\input amstex.tex
\documentstyle{amsppt}
\magnification=1200
\baselineskip=13pt
\hsize=6.5truein
\vsize=8.9truein
\countdef\sectionno=1
\countdef\eqnumber=10
\countdef\theoremno=11
\countdef\countrefno=12
\countdef\cntsubsecno=13
\sectionno=0
\def\newsection{\global\advance\sectionno by 1
                \global\eqnumber=1
                \global\theoremno=1
                \global\cntsubsecno=0
                {\bf\S}\the\sectionno.\ }

\def\newsubsection#1{\global\advance\cntsubsecno by 1
                     \xdef#1{{\S\the\sectionno.\the\cntsubsecno}}
                     \ \the\sectionno.\the\cntsubsecno.}

\def\theoremname#1{\the\sectionno.\the\theoremno
                   \xdef#1{{\the\sectionno.\the\theoremno}}
                   \global\advance\theoremno by 1}

\def\eqname#1{\the\sectionno.\the\eqnumber
              \xdef#1{{\the\sectionno.\the\eqnumber}}
              \global\advance\eqnumber by 1}

\global\countrefno=1

\def\refno#1{\xdef#1{{\the\countrefno}}\global\advance\countrefno by 1}

\def\thmref#1{#1}

\def\R{{\Bbb R}}

\def\C{{\Bbb C}}
\def\Z{{\Bbb Z}}
\def\Zp{{\Bbb Z}_+}
\def\Zt{{\Bbb Z}_2}
\def\Ztn{{\Bbb Z}_2^n}
\def\U{U_{q^{1/2}}({\frak{sl}}(2,\C))}
\def\a{\alpha}

\def\d{\delta}
\def\s{\sigma}
\def\t{\tau}

\def\e{\epsilon}
\def\vp{\varphi}
\def\p{\pi}
\def\r{\rho}

\refno{\ArikK}
\refno{\Bour}
\refno{\BrouCN}
\refno{\Cart}
\refno{\CharP}
\refno{\CurtIK}
\refno{\CurtRone}
\refno{\CurtRtwo}
\refno{\Dunk}
\refno{\DunkR}
\refno{\GaspR}
\refno{\GrozK}
\refno{\Hoef}
\refno{\Hump}
\refno{\Jimb}
\refno{\Koor}
\refno{\KoorZSE}
\refno{\Macd}
\refno{\MacdSB}
\refno{\Mats}
\refno{\StanAJM}
\refno{\StanGD}
\refno{\StanAKS}

\topmatter
\title $q$-Krawtchouk polynomials as spherical functions \\
on the Hecke algebra of type $B$
\endtitle
\rightheadtext{$q$-Krawtchouk polynomials and Hecke algebras}
\author H.T. Koelink\endauthor
\affil Universiteit van Amsterdam\endaffil
\address Vakgroep Wiskunde, Universiteit van Amsterdam,
Plantage Muidergracht 24,
1018 TV Amsterdam, the Netherlands\endaddress
\email koelink\@fwi.uva.nl\endemail
\date Report 96-07, May 22, 1996\enddate
\thanks Supported by the Netherlands
Organization for Scientific Research (NWO)
under project number 610.06.100. \endthanks
\abstract
The generic Hecke algebra for the hyperoctahedral group, i.e.
the Weyl group of type $B_n$, contains the generic Hecke
algebra for the symmetric group, i.e. the Weyl group
of type $A_{n-1}$, as a subalgebra. Inducing the index representation
of the subalgebra gives a Hecke algebra
module, which splits multiplicity free.
The corresponding zonal spherical functions
are calculated in terms of $q$-Krawtchouk polynomials.
The result covers a number of previously established
interpretations of ($q$-)Krawtchouk polynomials
on the hyperoctahedral group, finite groups
of Lie type, hypergroups and the quantum $SU(2)$ group.
\endabstract
\endtopmatter
\document


\subhead\newsection Introduction
\endsubhead

The hyperoctahedral group is the finite group of signed permutations, and
it contains the permutation group as a subgroup. The functions on the
hyperoctahedral group, which are left and right invariant with respect
to the permutation group, are spherical functions. The zonal
spherical functions are the spherical functions which are contained
in an irreducible subrepresentation of the group algebra under the
left regular representation. The zonal
spherical functions are known in terms of finite discrete orthogonal
polynomials, the so-called symmetric Krawtchouk polynomials. This
result goes back to Vere-Jones in the beginning of the seventies,
and is also contained in the work of Delsarte on coding theory,
see \cite{\Dunk} for information and references.

There are also group theoretic interpretations
of $q$-Krawtchouk polynomials. Here $q$-Krawtchouk polynomials
are $q$-analogues of the Krawtchouk polynomials in the sense that for
$q$ tending to one we recover the Krawtchouk polynomials. There
is more than one $q$-analogue for the Krawtchouk polynomial, but we
consider the $q$-analogue which is commonly called $q$-Krawtchouk
polynomial, see \S 7 
for its definition in terms of basic hypergeometric series \cite{\GaspR}.
In particular there is the interpretation by Stanton
\cite{\StanAJM}, \cite{\StanAKS}
of $q$-Krawtchouk polynomials (for specific
values of its parameter and of $q$) as spherical functions on
so-called finite groups of Lie type \cite{\Cart} which
have the hyperoctahedral group as Weyl group.

Outside the group theoretic setting there are some relevant
interpretations of the ($q$-)Krawtchouk polynomials.
The non-symmetric Krawtchouk polynomials have an interpretation
as symmetrised characters on certain hypergroups, cf.
Dunkl and Ramirez \cite{\DunkR}, and the hyperoctahedral group case is
recovered by a suitable specialisation. For the 
quantum $SU(2)$ group case the $q$-Krawtchouk polynomials enter
the picture as matrix elements of a basis transition in the
representation space of the irreducible
unitary representations, see Koornwinder \cite{\KoorZSE}.

The hyperoctahedral group is the Weyl group for the root system of type
$B$, and we can associate the (generic) Hecke algebra to it.
There is a subalgebra corresponding to the Hecke algebra for
the symmetric group, and we can consider functions which are
left and right invariant with respect to the index representation of
the Hecke algebra for the symmetric group. The index representation
is the analogue of the trivial representation.
We show that the zonal spherical functions on this Hecke
algebra can be expressed in terms of
$q$-Krawtchouk polynomials by deriving and solving a second-order
difference equation for the zonal spherical functions.
This interpretation of $q$-Krawtchouk polynomials gives a unified
approach to the interpretations of ($q$-)Krawtchouk polynomials as
zonal spherical functions on the hyperoctahedral group and the
appropriate finite groups of Lie type, since these can be
obtained by suitable specialisation. Moreover, it contains
the Dunkl and Ramirez result on the interpretation of
Krawtchouk polynomials as symmetrised characters on certain
hypergroups, and it gives a conceptual explanation
of the occurence of the $q$-Krawtchouk polynomials in the
quantum $SU(2)$ group setting.

The main reason why this procedure works is that the Hecke algebra
module induced from the index representation of the
subalgebra splits multiplicity free. This should be the
case, since the induced representation from the symmetric
group to the hyperoctahedral group splits multiplicity free.
To every Coxeter group a Hecke algebra can be associated
so we might consider each case of a finite Coxeter group
(here the hyperocathedral group)
with a maximal parabolic subgroup (here the
symmetric group) for which the induced representation from
the maximal parabolic subgroup splits multiplicity free.
A list of these situations can be found in
\cite{\BrouCN, Thm.~10.4.11}. In order to obtain
interesting spherical functions we have to restrict to
the cases $A_n$, $B_n=C_n$, $D_n$, for which
we have a parameter $n$, with $I_2^m$ as a possible
addition to this list, cf. \cite{\StanAKS, Table~1}.
For the Weyl group cases,
i.e. $A_n$, $B_n$, $D_n$, only the non-simply-laced $B_n$ is
interesting, since the corresponding Hecke algebra has two
parameters. The Hecke algebras for the simply-laced $A_n$ and $D_n$
only have one parameter, so we cannot expect an extension of the
results tabulated in Stanton \cite{\StanAKS, Table~2}.
The case studied here corresponds to $B_{n,n}$
(notation of \cite{\BrouCN, p.~295, Thm.~10.4.11}).

The contents of the paper are as follows. In \S 2 we
investigate the hyperoctahedral group in more detail.
In particular we study the coset representatives
of minimal length with respect to the permutation subgroup.
The Hecke algebra ${\Cal H}_n$ for the hyperoctahedral group and
the subalgebra ${\Cal F}_n$ corresponding to the Hecke algebra
for the symmetric group are introduced in \S 3.
The representation of ${\Cal H}_n$ obtained by inducing the
index representation of ${\Cal F}_n$ is also studied in
\S 3. The induced module $V_n$ is an analogue of $\C[\Ztn]$
and is a commutative algebra carrying a non-degenerate
bilinear form. In \S4 we let the quantised universal
enveloping algebra $\U$ for ${\frak{sl}}(2,\C)$ act on $V_n$,
which by Jimbo's analogue of the Frobenius-Schur-Weyl duality
gives the commutant of the action of ${\Cal F}_ n$ on $V_n$.
So the ${\Cal F}_n$-invariant elements, i.e. the elements
transforming according to the index representation of ${\Cal F}_n$,
in $V_n$ are identified with an irreducible module of $\U$.

Next, in \S 5 we calculate the characters of the algebra
$V_n$, and we give the corresponding orthogonality relations.
Using the non-degenerate bilinear form we identify
$V_n$ with its dual. The contragredient representation
is investigated in \S 6, and we decompose $V_n^\ast\cong V_n$
multiplicity free into
irreducible ${\Cal H}_n$-modules. In \S 6 we also investigate
the $\U$-module of ${\Cal F}_n$-invariant elements in greater
detail by giving the action of the generators of
$\U$ in an explicit basis.
A second-order difference equation as well as
orthogonality relations  for the zonal spherical
elements are derived in \S 7. The second-order
difference equation is used to identify the zonal
spherical element with $q$-Krawtchouk polynomials.
The relation with the interpretations of
($q$-)Krawtchouk polynomials alluded to in the
first few paragraphs is worked out in some more
detail in \S 7.

\demo{Notation}
The notation for $q$-shifted factorials
$(a;q)_k= \prod_{i=0}^{k-1} (1-aq^i)$,
$k\in \Zp$, and basic hypergeometric
series
$$
{}_{r+1}\varphi_r \left(
{{a_1,\ldots, a_{r+1}}\atop{b_1,\ldots, b_r}};q,z\right) =
\sum_{k=0}^\infty {{(a_1;q)_k \ldots (a_{r+1};q)_k}\over
{(b_1;q)_k \ldots (b_r;q)_k}} {{z^k}\over{(q;q)_k}}
$$
follows Gasper and Rahman \cite{\GaspR}. If one of the
upper parameters $a_i$ equals $q^{-d}$, $d\in\Zp$, then
the series terminates. If one of the lower parameters $b_j$
equals $q^{-n}$, $n\in\Zp$, then the series is not
well-defined unless one of the upper parameters
equals $q^{-d}$, $d\in\{ 0,1,\ldots,n\}$ with the
convention that if $d=n$ we consider it as a
terminating series of $n+1$ terms.
\enddemo

\demo{Acknowledgement} I thank Gert Heckman and Eric Opdam
for suggestions and their help in searching the literature.
\enddemo


\subhead\newsection The hyperoctahedral group
\endsubhead

The hyperoctahedral group is the semi-direct product
$H_n= \Ztn \rtimes S_n$, where
$\Zt=\{-1,1\}$ considered as the multiplicative group of
two elements and $S_n$ is the symmetric group, the group
of permutations on $n$ letters. So $H_n$ is the wreath
product $\Zt\wr S_n$.
The order of $H_n$ is $2^n\, n!$.

The hyperoctahedral group
is the Weyl group for the root system of type $B_n$ (or $C_n$), cf.
\cite{\Bour, Planche~II}. The hyperoctahedral group
operates on the standard $n$-dimensional Euclidean space
$\R^n$ with coordinates $\e_i$, $i=1,\ldots,n$.
The group $S_n$ acts by permutation of the coordinates and
$x\in\Ztn$ acts by sign changes, i.e. $x$ maps $\e_i$ to
$x_i\e_i$. Observe that $x\in\Ztn \subset H_n$ is involutive.
The roots, denoted by $R$,
are $\pm\e_i$, $1\leq i\leq n$, and $\pm\e_i\pm\e_j$,
$1\leq i<j\leq n$. The simple roots are
$\a_i=\e_i-\e_{i+1}$, $i=1,\ldots, n-1$, and $\a_n=\e_n$ and the
corresponding positive roots, denoted by $R^+$,
are $\e_i$, $1\leq i\leq n$, and $\e_i-\e_j$,
$1\leq i<j\leq n$, and $\e_i+\e_j$,
$1\leq i<j\leq n$. The other roots form the negative roots $R^-$.

Using this realisation of $H_n$ we see that for $x\in\Ztn$ and
$\s\in S_n$ we have $\s x\s^{-1} = x^\s\in \Ztn$, where
the action of $S_n$ on $\Ztn$ is given by
$$
x^\s = (x_{\s^{-1}(1)},\ldots,  x_{\s^{-1}(n)}).
\tag\eqname{\vglcommrelZnSn}
$$

The hyperoctahedral group is a Coxeter group \cite{\Bour}, \cite{\Hump}
with the generating set of reflections $S$ given by $s_i$,
$i=1,\ldots, n$, with $s_i$ the reflection in the positive
root $\a_i$. Thus $H_n$ is generated by $s_i$, $i=1,\ldots,n$,
subject to the relations $s_i^2=1$, and
$$
\gather
s_is_j=s_js_i, \ \ |i-j|>1,\qquad
s_is_{i+1}s_i=s_{i+1}s_is_{i+1},\ \ i=1,\ldots,n-2, \\
s_n s_{n-1}s_ns_{n-1} = s_{n-1}s_ns_{n-1}s_n.
\endgather
$$
In particular, the permutation group $S_n\subset H_n$
is the maximal parabolic subgroup generated by
$I=\{s_1,\dots,s_{n-1}\}\subset S$.
Hence there exist distinguished coset representatives for
$H_n/S_n$, cf. \cite{\Bour, p.~37, ex.~3}, \cite{\Hump, \S 1.10}.

\proclaim{Proposition \theoremname{\propdistcosetreprs}}
For $x\in\Ztn$ let $u_x\in H_n$ be the unique element of minimal
length in the coset $x\, S_n$. Then $u_x = x\, \s_x$
with
$$
\align
\ell(x) & = (1+2n)w(x) - 2\sum_{j=1}^{w(x)} i_j, \\
\ell\bigl(u_x\bigr) &= (1+n)w(x) - \sum_{j=1}^{w(x)} i_j, \\
\ell\bigl(\s_x\bigr) &= n\,w(x) - \sum_{j=1}^{w(x)} i_j,
\endalign
$$
where $w(x)= \#\{i\mid x_i=-1\}$ is the (Hamming) weight of
$x\in\Ztn$ and $x_{i_1}=\ldots=x_{i_{w(x)}}=-1$.
If $w\in xS_n$, then $w=u_x\s$ with
$\ell(w)=\ell\bigl(u_x\bigr) + \ell(\s)$.
\endproclaim

\demo{Proof} The last statement holds in greater generality, cf.
the general theorem on the distinguished coset representatives
\cite{\Bour, p.~37, ex.~3}, \cite{\Hump, \S 1.10} and it implies the
uniqueness of the coset representatives of minimal length.

Let $\s\in S_n$ be arbitrary, then
we can consider the action of $\s x$ on the negative roots and count
how many negative roots are made positive. Since for $w\in H_n$ we have
$\ell(w) = \# \bigl( R^+\cap w^{-1}R^-\bigr)$, we obtain
$$
\align
& \ell\bigl(x\s^{-1}\bigr) = \ell\bigl( x^{-1}\s^{-1}\bigr) =
\ell\bigl( (\s x)^{-1}\bigr) = \# \bigl( R^+\cap (\s x)R^-\bigr) =
 \# \{ i\mid x_i =-1\} + \\
& \# \{ i<j \mid \bigr( x_i=x_j=1 \land \s(i)>\s(j)\bigr) \, \lor\,
\bigl( x_i=-1 \land x_j=1\bigr) \,\lor\, \\
&\qquad\qquad\qquad \bigl( x_i=x_j=-1 \land \s(i)<\s(j)\bigr) \} \ + \\
& \# \{ i<j \mid \bigl( x_i=-1\land x_j=1\land \s(i)<\s(j)\bigr)\,\lor\,
\bigl( x_i=x_j=-1\bigr)\, \lor\, \\
&\qquad\qquad\qquad  \bigl( x_i=1\land x_j=-1\land
\s(i)>\s(j)\bigr) \},
\endalign
$$
where the three sets in the
last equality follow from counting the positive roots among respectively
$\s x(-\e_i)$, $\s x( -\e_i+\e_j)$, $\s x(-\e_i-e_j)$, $i<j$. We rewrite
this in a $\s$-dependent part and a $\s$-independent part as
$$
\aligned
 \ell\bigl(x\s^{-1}\bigr) = w(x) & + \#\{ i<j \mid x_i = -1\} \\
& + \# \{ i< j\mid x_i=x_j=1 \,\land\, \s(i)>\s(j)\} \\
& + \# \{ i< j\mid x_i=x_j=-1 \,\land\, \s(i)<\s(j)\} \\
& + \# \{ i< j\mid x_i=-1\,\land\, x_j=1 \,\land\, \s(i)<\s(j)\} \\
& + \# \{ i< j\mid x_i=1\,\land\, x_j=-1 \,\land\, \s(i)>\s(j)\} .
\endaligned
\tag\eqname{\vglexprlengthxs}
$$
In particular, if $\s=1$ we obtain
$$
\ell(x) = w(x) + 2\#\{ i<j\mid x_i=-1\} = w(x)
+ 2\sum_{j=1}^{w(x)} \bigl( n-i_j\bigr).
$$

The $\s$-dependent part in \thetag{\vglexprlengthxs} is zero, and
thus minimal, for $\s$ defined inductively by $\s(1)=n$ if $x_1=-1$
and $\s(1)=1$ if $x_1=1$, and $\s(i)$ is as large, respectively small,
as possible if $x_i=-1$, respectively $x_i=1$. Now define $\s_x$ to
be the inverse of the $\s$ defined in this way, and let $u_x=x\, \s_x$.
Then $u_x=x\, \s_x$ has length $\ell\bigl(u_x\bigr)
= w(x) + \#\{ i<j\mid x_i=-1\}$, which is minimal in the coset
$x\, S_n$ by \thetag{\vglexprlengthxs}. The statement on
$\ell\bigl(\s_x\bigr)$ follows directly, or from the general last
statement of the proposition.
\qed\enddemo

\demo{Remark \theoremname{\rempropdistcosetreprs}} (i) If we write $u_x$
as the $n\times n$ signed permutation matrix, then it is characterised
by the conditions (i) all $-1$'s occur columnwise to the right of all
$+1$'s, (ii) the $-1$'s decrease in columns as they increase by rows,
(iii) the $+1$'s increase in columns as they increase by rows, (iv)
the non-zero entry in the $i$-th row is $x_i$. The
permutation matrix for $\s_x$ is obtained from the one for $u_x$
by replacing all $-1$'s by $+1$'s.

(ii) The cardinality of $H_n/S_n$ is $2^n$ and Proposition
\thmref{\propdistcosetreprs} shows how this can be pa\-ra\-me\-tri\-sed
by $\Ztn$.

(iii) Proposition \thmref{\propdistcosetreprs} gives the length of
$x\in\Ztn$. To give an explicit reduced expression we first observe
that $x=x^{i_1}\ldots x^{i_{w(x)}}$ is a decomposition in
commuting elements, where $x^j$ is defined by $w(x^j)=1$ and
$i_1(x^j)=j$, or $x^j=(1,\ldots,1,-1,1,\ldots,1)$ with $-1$ at the
$j$-th spot. Proposition \thmref{\propdistcosetreprs} implies that
$\ell(x) = \sum_{j=1}^{w(x)} \ell(x^{i_j})$, so that a reduced
expression for $x$ is obtained by inserting a reduced expression
for each $x^j$. Finally, note that
$$
x^j = s_js_{j+1}\ldots s_{n-1}s_ns_{n-1}s_{n-2}\ldots s_j
\tag\eqname{\vglreducedexprxj}
$$
is a reduced expression.
\enddemo

\proclaim{Proposition \theoremname{\propcommrelux}} Let $u_x$,
$x\in\Ztn$, be the coset representatives of minimal length as in
Proposition \thmref{\propdistcosetreprs}, then
$$
\gather
\ell\bigr( s_i u_x\bigr) = \ell\bigl( u_x\bigr) + 1
\Longleftrightarrow x_i \geq x_{i+1}, \quad 1\leq i<n, \\
\ell\bigr( s_n u_x\bigr) = \ell\bigl( u_x\bigr) + 1
\Longleftrightarrow x_n=1.
\endgather
$$
Moreover, for $1\leq i<n$ we have $s_iu_x = u_{x^{s_i}}$ if
$x_i\not=x_{i+1}$ and $s_i u_x = u_xs_j$ if $x_i=x_{i+1}$ with
$j=\min \bigl( \s_x^{-1}(i), \s_x^{-1}(i+1)\bigr)$ and
$s_n u_x = u_{xx^n}$.
\endproclaim

\demo{Proof} Recall, see \cite{\Hump, \S\S 1.6-7}, that
$\ell(s_iw)=\ell(w)+1$ if and only of $w^{-1}\a_i\in R^+$. So,
taking $w=u_x$, we have to calculate for $1\leq i<n$,
$$
u_x^{-1}\,\a_i = \s_x^{-1} x\, \a_i = x_i \e_{\s_x^{-1}(i)} - x_{i+1}
\e_{\s_x^{-1}(i+1)}.
$$
For $x_i=-1$, $x_{i+1}=1$ this root is negative, and for
$x_i=1$, $x_{i+1}=-1$ this root is positive. By construction of the
permutation $\s_x$, cf. the proof of Proposition
\thmref{\propdistcosetreprs}, we have $\s_x^{-1}(i)<\s_x^{-1}(i+1)$ if
$x_i=x_{i+1}=1$ and $\s_x^{-1}(i)>\s_x^{-1}(i+1)$ if
$x_i=x_{i+1}=-1$, so that in these cases $u_x^{-1}\,\a_i$ is a positive
root as well. Similarly,
$$
u_x^{-1}\, \a_n = \s_x^{-1} x\, \e_n = x_n \, \e_{\s_x^{-1}(n)}
$$
and this is a positive root if and only if $x_n=1$.

To prove the second statement recall, cf. \thetag{\vglcommrelZnSn},
that $\s x =x^\s\s$ for $x\in\Ztn$, $\s\in S_n$. Hence, for $1\leq i<n$
we have
$$
s_iu_x=s_i x\s_x = x^{s_i}s_i\s_x =
u_{x^{s_i}} \s_{x^{s_i}}^{-1} s_i\s_x
\in u_{x^{s_i}} S_n.
$$
Now we have to consider some cases separately. Firstly, if $x_i=x_{i+1}$,
then $x^{s_i}=x$, and by the part of the proposition already proved,
$$
\ell\big(s_iu_x\bigr) = \ell\bigl(u_x\bigr) + 1 =
\ell\bigl(u_{x^{s_i}}\bigr) + 1.
$$
Hence by Proposition \thmref{\propdistcosetreprs}
$\ell\bigl(\s_x^{-1} s_i\s_x\bigr)=1$, so that $\s_x^{-1} s_i\s_x=s_j$
for some $1\leq j\leq n$. Since $\s_x^{-1} s_i\s_x$ interchanges
$\e_{\s_x^{-1}(i)}$ and $\e_{\s_x^{-1}(i+1)}$ and fixes the other
coordinates $\e_k$ and since
$\s_x^{-1}(i)$ and $\s_x^{-1}(i+1)$ differ by $1$ by construction of
$\s_x$, it follows that $j=\s_x^{-1}(i)$ if $x_i=x_{i+1}=1$ and
$j=\s_x^{-1}(i+1)$ if $x_i=x_{i+1}=-1$ by construction of
$\s_x$, or $j=\min \bigl(
\s_x^{-1}(i), \s_x^{-1}(i+1)\bigr)$.

In case $x_i=1$, $x_{i+1}=-1$ we have by Proposition
\thmref{\propdistcosetreprs} and by the part of the proposition
already proved
$$
\ell\bigl( u_{x^{s_i}}\bigr) = \ell\bigl(u_x\bigr) + 1 =
\ell\bigl(s_i u_x\bigr).
$$
Since $u(x^{s_i})$ and $s_i u_x$ are in the same coset $x^{s_i}S_n$
and have the same length, it follows that $s_i u_x=u(x^{s_i})$.
Replace $x$ by $x^{s_i}$ to find the same conclusion in case
$x_i=-1$, $x_{i+1}=1$.

It remains to consider
$$
s_n u_x = s_n x \s_x = xx^n\s_x = u_{xx^n} \s_{xx^n}^{-1}\s_x
\in u_{xx^n}S_n.
$$
By Proposition \thmref{\propdistcosetreprs}
and by the part of the proposition already proved we see that for
both $x_n=1$ and $x_n=-1$ we have
$\ell\bigl(u_{xx^n}\bigr) = \ell\bigl(s_nu_x\bigr)$
and thus $s_n u_x=u_{xx^n}$ by uniqueness.
\qed\enddemo


\subhead\newsection The Hecke algebra for $H_n$ and an induced
representation
\endsubhead

The (generic) Hecke algebra can be associated to any Coxeter group,
cf. \cite{\Bour, p.~54, ex.~22}, \cite{\CurtRtwo,
\S 68A}, \cite{\Hump, \S 7.1}. In particular, we define the Hecke
algebra ${\Cal H}_n= {\Cal H}(H_n)$ associated with the
hyperoctahedral group as the algebra over the field
$\C(p^{1/2},q^{1/2})$ with elements $T_w$, $w\in H_n$,
subject to the relations
$$
\aligned
&T_i T_w = T_{s_iw}, \qquad
\text{if}\quad \ell(s_iw)>\ell(w), \\
&T_i^2 = (q-1)T_i + q, \qquad 1\leq i <n, \\
&T_n^2 = (p-1)T_n + p,
\endaligned
\tag\eqname{\vgldefHeckeBn}
$$
where we use the notation $T_i=T_{s_i}$, $1\leq i\leq n$. The identity
of ${\Cal H}_n$ is $1=T_e$ with $e\in H_n$ the identity of the
hyperoctahedral group. Note that the elements $T_i$, $1\leq i\leq n$,
generate ${\Cal H}_n$.

The Hecke algebra ${\Cal H}_n$ can be defined over the ring $\Z[p,q]$,
but then we would have to extend it later.

Similarly, we define the Hecke algebra
$\tilde {\Cal F}_n ={\Cal H}(S_n)$ as the algebra over $\C(q^{1/2})$
with generators $T_\s$, $\s\in S_n$, subject to the relations
$$
\align
&T_i T_\s = T_{s_i\s}, \qquad \text{if}\quad \ell(s_i\s)>\ell(\s),
\quad 1\leq i <n, \\
&T_i^2 = (q-1)T_i + q, \qquad 1\leq i <n,
\endalign
$$
so that we may view
${\Cal F}_n = \C(p^{1/2},q^{1/2})\otimes_{\C(q^{1/2})} \tilde{\Cal F}_n$
as a (maximal parabolic) subalgebra of ${\Cal H}_n$.

The Hecke algebra ${\Cal F}_n$ has two one-dimensional
representations, cf. \cite{\CurtIK, \S 10}, namely the index
representation
$$
\iota\colon {\Cal F}_n \to \C(p^{1/2},q^{1/2}),
\quad \iota(T_i)=q,\ 1\leq i<n,
$$
and the sign representation $T_w\mapsto (-1)^{\ell(w)}$.
Observe that $\iota(T_\s)=q^{\ell(\s)}$.

The Hecke algebra ${\Cal H}_n$ has four one-dimensional representations
of which two representations restricted to ${\Cal F}_n$ give the
index representation $\iota$ of ${\Cal F}_n$. They are defined by
$$
\iota,\iota^\prime \colon {\Cal H}_n \to \C(p^{1/2},q^{1/2}),
\quad \iota|_{{\Cal F}_n}
=\iota,\ \iota^\prime|_{{\Cal F}_n}=\iota,\quad \iota(T_n)=p,\
\iota^\prime(T_n)=-1.
$$
Again, $\iota$ is called the index representation of ${\Cal H}_n$ and
denoted by the same symbol. The other two one-dimensional
representations of ${\Cal H}_n$ when restricted to ${\Cal F}_n$
give the sign representation of ${\Cal F}_n$, see
\cite{\CurtIK, \S 10}. The complete representation theory of
${\Cal H}_n$ and ${\Cal F}_n$ can be found in
Hoefsmit's thesis \cite{\Hoef}.

By $V_n$ we denote the induced module obtained from inducing the
one-dimensional representation $\iota$ of ${\Cal F}_n$ to ${\Cal H}_n$.
So
$$
V_n = {\Cal H}_n \otimes_{{\Cal F}_n} \C(q^{1/2}) =
\text{Ind}_{{\Cal F}_n}^{{\Cal H}_n}\ \iota,
$$
and ${\Cal H}_n$ acts from the left by multiplication. We denote the
corresponding representation of ${\Cal H}_n$ in $V_n$ by $\r$.
Define the mapping
$$
\p \colon {\Cal H}_n \to V_n,\qquad \p(T) = \r(T)(1\otimes 1) =
T\otimes 1.
$$

\proclaim{Theorem \theoremname{\propalgebrastructureV}}
{\rm (i)}
$\{ u(x)=\p(T_{u_x}) \mid x\in\Ztn\}$ forms a basis for $V_n$. With
respect to this basis the action of the generators of ${\Cal H}_n$
is given by, $1\leq i <n$,
$$
\align
\r(T_i)\, u(x) &=
\cases u(x^{s_i}), &\text{if $x_i=1$, $x_{i+1}=-1$,}\\
q\, u(x), &\text{if $x_i=x_{i+1}$,}\\
(q-1)\,u(x) + q\, u(x^{s_i}), &\text{if $x_i=-1$, $x_{i+1}=1$,}
\endcases \\
\r(T_n)\, u(x) &=
\cases u(xx^n), &\text{if $x_n=1$,}\\
(p-1)\, u(x) + p\, u(xx^n), &\text{if $x_n=-1$.}
\endcases
\endalign
$$

\noindent
{\rm (ii)} Define
$$
u(x) u(y) = q^{-\ell(\s_x) - \ell(\s_y)}
\p(T_xT_y),
\tag\eqname{\vgldefalgstruct}
$$
then $V_n$ is a commutative algebra over $\C(p^{1/2},q^{1/2})$.
As an algebra $V_n$ is generated
by $u(x^j)$, $1\leq j\leq n$, with $x^j$ defined in Remark
\thmref{\rempropdistcosetreprs}(iii), subject to the relations
$$
\aligned
u(x^j) u(x^k) &= u(x^k)u(x^j) = u(x^kx^j), \qquad k\not= j,\\
\bigl( u(x^j)\bigr)^2 &=
(q-1) u(x^j)\Bigl( u(x^{j+1}) + \ldots + u(x^n)\Bigr) +
(p-1) u(x^j) + pq^{n-j} u(1).
\endaligned
\tag\eqname{\vglcommrelV}
$$

\noindent
{\rm (iii)} The $\C(p^{1/2},q^{1/2})$-linear form $\t$ on $V_n$
defined by
$$
\t\bigl( u(x)\bigr) = \d_{x,1} = \d_{x_1,1}\ldots\d_{x_n,1}
$$
induces a non-degenerate symmetric associative bilinear form
$B(\xi,\eta) = \t(\xi\eta)$ on $V_n$ and
$$
B\bigl(u(x), u(y)\bigr) = \d_{x,y}\,  \iota(T_{u_x}).
$$
So $\{u(x)\}$ and $\{v(x)= \iota(T_{u_x})^{-1}u(x)\}$ are dual
bases for the bilinear form $B$.
\endproclaim

\demo{Remark \theoremname{\rempropalgebrastructureV}} (i)
The product in $V_n\cong {\Cal H}_nP$, with $P$ the central
idempotent of \thetag{4.6}, 
is defined by $(TP)(SP)=(TS)P$ instead of $(TP)(SP)=(TPS)P$.

(ii) Instead of the basis $\pi(T_{u_x})$ we can also use the basis
$\pi(T_x)=q^{\ell(\s_x)}\pi(T_{u_x})$, by Proposition
\thmref{\propdistcosetreprs}, for the induced module $V_n$.
We obtain the same action of the generators for ${\Cal H}_n$
as in Theorem \thmref{\propalgebrastructureV}(i) by using
Lemma~3.3 of Ariki and Koike \cite{\ArikK}.
\enddemo

\demo{Proof} (i) Since ${\Cal H}_n\cong \C[H_n]$ and
${\Cal F}_n\cong \C[S_n]$ and
the index representation corresponds to the trivial representation,
it follows that $V_n\cong \C[\Ztn] = \text{Ind}_{S_n}^{H_n} 1$. Hence,
$V_n$ is $2^n$-dimensional and Proposition
\thmref{\propdistcosetreprs} implies that
$\{ u(x) \mid x\in\Ztn\}$ forms a basis for $V_n$. The action of the
generators $T_i$, $1\leq i\leq n$, follow directly from the definition
of $V_n$, Proposition \thmref{\propcommrelux} and the relations
$$
\align
T_iT_w &= (q-1) T_w + qT_{s_iw},\
\text{if}\ \ell(s_iw) < \ell(w), \
1\leq i < n, \\
T_nT_w &= (p-1) T_w + pT_{s_nw},\
\text{if}\ \ell(s_nw) < \ell(w)
\endalign
$$
in ${\Cal H}_n$.

(ii) Now \thetag{\vgldefalgstruct} is
equivalent to $\p(T_x)\p(T_y)=\p(T_xT_y)$ and it
defines an algebra structure on $V_n$, since
${\Cal H}_n$ is an algebra and $\p$ is an algebra morphism.
To show
that it is a commutative algebra we remark that $T_xT_y=T_yT_x$
for $x,y\in\Ztn$. Since
$$
T_x = T_{x^{i_1}}\ldots T_{x^{i_{w(x)}}},
$$
see Remark \thmref{\rempropdistcosetreprs}(iii), this follows from
$T_{x^j}T_{x^k}=T_{x^k}T_{x^j}$. This is obvious for $j=k$. For
$j\not= k$ we have from Proposition \thmref{\propdistcosetreprs}
$\ell(x^j)+\ell(x^k) = \ell(x^jx^k)$ implying
$T_{x^j}T_{x^k}=T_{x^jx^k}=T_{x^kx^j}=T_{x^k}T_{x^j}$, since $\Ztn$ is
commutative.

To see that $V_n$ is generated by $u(x^j)$, $1\leq j\leq n$, we observe
that
$$
u(x^{i_1})\ldots u(x^{i_{w(x)}}) = q^{-\ell(\s_x)}
\p(T_x),
$$
which is proved by induction on the distance of $x$ and Proposition
\thmref{\propdistcosetreprs}. By the same Proposition we have
$T_x=T_{u_x} T_{\s_x^{-1}}$, and hence
$\p(T_x) = q^{\ell(\s_x)}u(x)$. So we conclude that
$$
u(x^{i_1})\ldots u(x^{i_{w(x)}}) = u(x).
\tag\eqname{\vglontbindingbasis}
$$
In order to complete the algebra structure of $V_n$
it remains to calculate
$$
\bigl(u(x^j)\bigr)^2 = q^{-2(n-j)} \p(T_{x^j}T_{x^j}).
$$
For $j=n$ we obtain the result from $T_n^2=(p-1)T_n+p$. In general we
use the reduced expression \thetag{\vglreducedexprxj} and
$T_j^2=(q-1)T_j+q$ to find
$$
T_{x^j}T_{x^j} = (q-1)T_{x^j}T_{x^{j+1}}T_j +
q T_j T_{x^{j+1}}T_{x^{j+1}} T_j.
$$
Hence, taking into account Proposition
\thmref{\propdistcosetreprs}, we find
$$
\bigl(u(x^j)\bigr)^2 = (q-1) u(x^j) u(x^{j+1})
+ \r(T_j)\, \bigl(u(x^{j+1})\bigr)^2.
$$
Since the case $j=n$ is true, we can use downward induction on
$j$ to find
$$
\multline
\bigl(u(x^j)\bigr)^2 = (q-1) u(x^j) u(x^{j+1})
+ (q-1) \sum_{p=j+2}^n \r(T_j) u(x^{j+1}x^p) \\
+ (p-1) \r(T_j) u(x^{j+1}) + pq^{n-j-1} \r(T_j) u(1)
\endmultline
$$
by \thetag{\vglontbindingbasis}. Now use part (i) of the
theorem to find \thetag{\vglcommrelV}.

(iii) The bilinear form is associative by construction and
symmetric since $V_n$ is a commutative algebra. Since the
last equality implies that $B$ is non-degenerate, it suffices
to prove
$$
\t \bigl( u(x)u(y)\bigr) = \d_{x,y} \iota(T_{u_x}).
$$

This can be proved by an induction argument. First we use
induction on $w(y)$, the case $w(y)=0$ being trivial. Let
$k=i_{w(y)}$ so that $y_k=-1$ and $y_{k+1}=\ldots = y_n=1$.
There are two cases to be considered, namely $x_k=1$ or $x_k=-1$.
If $x_k=1$ we have
$$
\t \bigl( u(x)u(y)\bigr) =
\t \bigl( u(xx^k)u(yx^k))\bigr) =
0 =\d_{x,y} \iota(T_{u_x})
$$
by \thetag{\vglontbindingbasis} and $w(yx^k)=w(y)-1$.

If $x_k=-1$ we use downward induction on $k$. So let us first consider
$k=n$. Then
$$
\t \bigl( u(x) u(y)\bigr) =
\t \bigl( u(xx^n)u(yx^n) u(x^n)^2 \bigr) =
(p-1)\, \t \bigl( u(x)u(yx^n) \bigr)
+ p\, \t \bigl( u(xx^n)u(yx^n)\bigr)
$$
by \thetag{\vglontbindingbasis} and \thetag{\vglcommrelV}.
By the induction hypothesis on $w(y)$, the first term equals zero
and the second term equals $p\d_{xx^n,yx^n} \iota(T_{u_{xx^n}}) =
\d_{x,y} \iota(T_{u_x})$ by Proposition \thmref{\propdistcosetreprs}.

Now we have
$$
\multline
\t \bigl( u(x)u(y)\bigr) =
\t \bigl( u(xx^k)u(yx^k)
u(x^k)^2 \bigr) = \\
(p-1)\, \t \bigl( u(x)u(yx^k)\bigr) +
\sum_{l=k+1}^n \t \bigl( u(x)u(yx^kx^l)\bigr)
+pq^{n-k}\, \t \bigl( u(xx^k) u(yx^k)\bigr).
\endmultline
$$
The first term is zero by the induction hypothesis for
$w(y)$ and all terms in the sum
are zero by the induction hypothesis on $k$ and the case $x_k=1$
already proved. By the induction hypothesis on $w(y)$
we obtain $pq^{n-k} \d_{xx^k,yx^k} \iota(T_{u_{xx^k}}) =
\d_{x,y} \iota(T_{u_x})$ by Proposition \thmref{\propdistcosetreprs}.
\qed\enddemo

The last statement of Theorem \thmref{\propalgebrastructureV}
is the analogue of the more general statement that
$T_w$ and $\iota(T_w)^{-1}T_{w^{-1}}$ are dual bases for the bilinear
form associated to the linear form $T_w\mapsto \d_{w,1}$. This holds
for the Hecke algebra associated to any finite Coxeter group.


\subhead\newsection ${\Cal F}_n$-invariant elements in the induced
representation
\endsubhead

Let us first introduce the orthonormal basis of $V_n$ defined by
$$
\hat u(x) = u(x) \bigl(\iota(T_{u_x})\bigr)^{-1/2},
\qquad x\in\Ztn,
$$
so that $B\bigl( \hat u(x),\hat u(y)\bigr)=\d_{x,y}$ or $\hat u(x)$,
$x\in\Ztn$, is an orthonormal basis of $V_n$ with respect to the
bilinear form $B$, which defines an inner product on $V_n$.
Then we have, cf. Theorem \thmref{\propalgebrastructureV}(i),
$$
\r(T_i)\, \hat u(x)=
\cases q^{1/2}\, \hat u(x^{s_i}), &\text{if $x_i=1$,
$x_{i+1}=-1$,}\\
q \,\hat u(x), &\text{if $x_i=x_{i+1}$,}\\
(q-1)\, \hat u(x) + q^{1/2}\, \hat u(x^{s_i}),
&\text{if $x_i=-1$, $x_{i+1}=1$.}
\endcases
\tag\eqname{\vglTiopvx}
$$
Let $y\in\Zt^{i-1}$, $z\in\Zt^{n-i-1}$ fixed and define
the ordered basis $f_{-1}\otimes f_{-1}=\hat u(y,-1,-1,z)$,
$f_{-1}\otimes f_{1}=\hat u(y,-1,1,z)$,
$f_{1}\otimes f_{-1}=\hat u(y,1,-1,z)$,
$f_{1}\otimes f_{1}=\hat u(y,1,1,z)$. Then $T_i$ leaves the space
spanned by these four elements invariant
and with respect to this basis $\r(T_i)$ is
represented by the $4\times 4$-matrix
$$
R = \pmatrix q & 0 & 0 & 0 \\
0 & q-1 & q^{1/2} & 0 \\
0 & q^{1/2} & 0 & 0 \\
0 & 0 & 0 & q \endpmatrix,
\tag\eqname{\vglRmatrixTi}
$$
which is closely related to the $R$-matrix in the fundamental
representation for the quantised universal enveloping algebra
$\U$. So let us recall the definition of $\U$, see e.g.
\cite{\CharP, Def.~9.1.1}.

\proclaim{Definition \theoremname{\defquantumsl2}}
$\U$ is the associative algebra with unit over $\C(q^{1/2})$
with generators $K$, $K^{-1}$, $E$ and $F$ subject to the relations
$$
\gather
KK^{-1}=1=K^{-1}K, \quad KE=qEK,\quad KF=q^{-1}FK,\\
EF-FE = {{K-K^{-1}}\over{q^{1/2}-q^{-1/2}}}.
\endgather
$$
There exists a Hopf-algebra structure on $\U$ with the
comultiplication
$$
\Delta \colon \U \to \U\otimes \U,
$$
which is an algebra homomorphism, given by
$$
\gather
\Delta(K) = K\otimes K, \qquad \Delta(E) =
K\otimes E + E\otimes 1,\\
\Delta(F) = 1\otimes F + F\otimes K^{-1},\qquad
\Delta(K^{-1}) = K^{-1}\otimes K^{-1}.
\endgather
$$
\endproclaim

\demo{Remark \theoremname{\remdefquantumsl2}} The correspondence
with the quantised universal enveloping algebra as in
\cite{\KoorZSE, \S 3} is obtained by identifying $K$,
$K^{-1}$, $E$ and $F$ with $A^2$, $D^2$, $q^{-1/4}AB$,
$q^{1/4}CD$ of \cite{\KoorZSE, \S 3}, where we have replaced
$q$ by $q^{1/2}$.
\enddemo

The representation theory of $\U$ is well-known, cf. e.g.
\cite{\CharP, \S 10.1}, \cite{\KoorZSE, \S3},
and is similar to the representation theory
of the Lie algebra ${\frak{sl}}(2,\C)$. We recall some of these
results in the next theorem.

\proclaim{Theorem \theoremname{\thmreprquantumsl2}}
{\rm (i)} There is precisely one irreducible
$\U$-module $W_N$ of each dimension $N+1$ over $\C(q^{1/2})$
with highest weight vector $v_+$, i.e. $K\cdot v_+ = q^{N/2} v_+$,
$E\cdot v_+ = 0$.

\noindent
{\rm (ii)} The Clebsch-Gordan decomposition holds; as $\U$-modules
$$
W_N\otimes W_M = \bigoplus_{k=0}^{\min(N,M)} W_{M+N-2k}.
$$
\endproclaim

The tensor product representation of $\U$ is defined using the
comultiplication $\Delta$. If $\nu^1$, $\nu^2$ are two
representations of $\U$ acting in $W^1$, $W^2$, then the
representation $\nu^1\otimes \nu^2$
acts in $W^1\otimes W^2$ and is defined by
$$
\bigl(\nu^1\otimes \nu^2\bigr)(X)\, w^1\otimes w^2 =
\sum_{(X)} \nu^1(X_{(1)})\,w^1\otimes  \nu^2(X_{(2)})\, w^2,
$$
where $\Delta(X) = \sum_{(X)} X_{(1)}\otimes X_{(2)}$.

We can define an antilinear $\ast$-operator on $\U$ by
$(q^{1/2})^\ast=q^{1/2}$ and
$$
K^\ast=K, \quad E^\ast=q^{-1/2}FK,\quad F^\ast =q^{1/2}K^{-1}E,
\quad (K^{-1})^\ast=K^{-1}.
\tag\eqname{\vglatstU}
$$
There are other $\ast$-structures possible on $\U$ corresponding to
different real forms of ${\frak{sl}}(2,\C)$. This $\ast$-structure
corresponds to the compact real form ${\frak{su}}(2)$ and
for this $\ast$-structure the irreducible $\U$-modules described
in Theorem \thmref{\thmreprquantumsl2}, and all the
tensor product representations built from the irreducible
representations, are $\ast$-representations
of $\U$. For more information on this subject we refer to
\cite{\CharP, Ch.~9}.

Let us now consider the fundamental 2-dimensional representation
in $W=W_1$ of $\U$. Let $\{e_{-1},e_1\}$ be the standard orthonormal
basis of $W$, then we have
$$
K \mapsto \pmatrix q^{1/2} & 0\\ 0 & q^{-1/2}\endpmatrix, \quad
E \mapsto \pmatrix 0 & 1 \\ 0 & 0\endpmatrix, \quad
F \mapsto\pmatrix 0 & 0 \\ 1 & 0 \endpmatrix.
\tag\eqname{\vglactiefundrep}
$$
Consider the $n$-fold tensor product representation $t$
in $W^{\otimes n}$ and extend $t$ to
a representation of $\U$ in
$\C(p^{1/2},q^{1/2})\otimes_{\C(q^{1/2})}W^{\otimes n}$ by
$t(X)\bigl(f\otimes w\bigr) = f\otimes t(X)w$.
By identifying
$\C(p^{1/2},q^{1/2})\otimes_{\C(q^{1/2})}W^{\otimes n}$
with $V_n$ by the unitary mapping
$$
f\otimes e_{i_1}\otimes e_{i_2} \otimes \ldots
\otimes e_{i_n} \longmapsto
f\,  \hat u (i_1,i_2,\ldots, i_n), \qquad f\in \C(p^{1/2},q^{1/2}),
\tag\eqname{\vglidenWnenVn}
$$
we have the representation $t$ of $\U$ acting in $V_n$.
This representation is a $\ast$-representation of $\U$ for the
$\ast$-operator defined by \thetag{\vglatstU}.
It is straightforward to check that for $n=2$ the representation
$t$ of $\U$ in $V_2$ commutes with
the matrix $R$ of \thetag{\vglRmatrixTi} and this is part of Jimbo's
theorem on the analogue of the Frobenius-Schur-Weyl duality, cf.
\cite{\Jimb, Prop.~3}, \cite{\CharP, \S 10.2.B}.

\proclaim{Theorem \theoremname{\thmJimbo}}{\rm (Jimbo)}
The algebras $t\bigl(\U\bigr)$ and $\r({\Cal F}_n)$
are each others commutant in $\hbox{\rm End}_{\C(p^{1/2},q^{1/2})}(V_n)$.
\endproclaim

Jimbo's theorem can be used to determine the ${\Cal F}_n$-invariant
elements in $V_n$. We call an element $v\in V_n$ a
${\Cal F}_n$-invariant element if
$\r(T_\s)\, v = \iota(T_\s)\, v = q^{\ell(\s)}\,v$
for all $\s\in S_n$. So $v$ is ${\Cal F}_n$-invariant if it
realises the one-dimensional index representation of ${\Cal F}_n$.
Let us define the corresponding idempotent, central in ${\Cal F}_n$,
$$
P= {1\over{P_A(q)}} \sum_{\s\in S_n} T_\s
\in {\Cal F}_n\subset {\Cal H}_n.
\tag\eqname{\vgldefP}
$$
Then in any representation of ${\Cal H}_n$ the operator
corresponding to $P$ acts as a projection operator
on the ${\Cal F}_n$-invariant elements since
$P(T_i-q)=(T_i-q)P=0$, $1\leq i<n$. Here
$P_A(q)=\sum_{\s\in S_n} q^{\ell(\s)}$
is the Poincar\'e polynomial for the Coxeter group $S_n$.
So for arbitrary $v\in V_n$ the Hecke symmetrised vector
$$
\r(P)\, v = {1\over{P_A(q)}} \sum_{\s\in S_n} \r(T_\s)\, v
$$
is ${\Cal F}_n$-symmetric.

\proclaim{Proposition \theoremname{\propFinvariantelts}}
The space of ${\Cal F}_n$-invariant elements in $V_n$
realises the unique irreducible $n+1$-dimensional $\U$-module
with highest weight vector $u(-1,-1,\ldots,-1)$.
\endproclaim

\demo{Proof}
By Jimbo's Theorem \thmref{\thmJimbo} we see that the space
of ${\Cal F}_n$-invariant elements is an invariant $\U$-module.
Since $u(-1,-1,\ldots,-1)$ is a ${\Cal F}_n$-invariant vector by
\thetag{\vglTiopvx}, it suffices to check
$$
\aligned
t(K)\, u(-1,-1,\ldots,-1) &= q^{n/2}\, u(-1,-1,\ldots,-1), \\
t(E)\, u(-1,-1,\ldots,-1) &= 0,
\endaligned
\tag\eqname{\vglABactietebewijzen}
$$
by Theorem \thmref{\thmreprquantumsl2}(i). Let
$\Delta^{(2)} = \Delta$ and inductively
$$
\Delta^{(n)} = (1\otimes \ldots\otimes
1\otimes \Delta)\circ \Delta^{(n-1)}
\colon \U \to \Bigl( \U\Bigr)^{\otimes n},
$$
then it follows by induction on $n$ from
Definition \thmref{\defquantumsl2} that
$$
\aligned
\Delta^{(n)}(K) &= K\otimes K\otimes K
\otimes \ldots \otimes K, \\
\Delta^{(n)}(E) &= \sum_{i=1}^n K\otimes
\ldots \otimes K\otimes E \otimes 1
\otimes \ldots \otimes 1,
\endaligned
\tag\eqname{\vglDeltanopAB}
$$
where the $E$ occurs at the $i$-th component 
of the tensor product. From \thetag{\vglactiefundrep} 
and the identification
\thetag{\vglidenWnenVn}
we directly obtain \thetag{\vglABactietebewijzen}, since
$K\cdot e_{-1}=q^{1/2}e_{-1}$ and $E\cdot e_{-1}=0$.

{}From iteration of the Clebsch-Gordan decomposition in
Theorem \thmref{\thmreprquantumsl2} it follows that the
$n+1$-dimensional irreducible module $W_n$ occurs with
multiplicity one in
$W^{\otimes n}$, which implies uniqueness.
\qed\enddemo

In \S 6 
we give an explicit description of this
$\U$-module of ${\Cal F}_n$-invariant
elements.


\subhead\newsection Characters of $V_n$ and orthogonality relations
\endsubhead

Since $V_n$ is a commutative algebra, by Theorem
\thmref{\propalgebrastructureV}(ii), all irreducible
representations are one-dimensional, so we now investigate the
characters of $V_n$. There are $2^n$ characters and they span
$V_n^\ast=\text{Hom}_{\C(p^{1/2},q^{1/2})}(V_n, \C(p^{1/2},q^{1/2}))$.

\proclaim{Theorem \theoremname{\thmcharactersV}}
$V_n^\ast$ is spanned by the characters $\chi_y$, $y\in\Ztn$,
defined by
$$
\chi_y\bigl(u(x^j)\bigr) =
\cases pq^{m_j(y)}, &\text{if $y_j=1$,} \\
-q^{n-j-m_j(y)}, &\text{if $y_j=-1$,}
\endcases
$$
where $m_j(y)=\#\{ p>j\mid y_p=1\}$.
\endproclaim

\demo{Proof} The $\chi_y$ are obviously different for different
$y\in\Ztn$. It suffices to check that $\chi_y$ preserves the
quadratic relations of \thetag{\vglcommrelV}. So we calculate
$$
\sum_{k=j+1}^n \chi_y\bigl(u(x^k)\bigr) =
\mathop{\sum_{k=j+1}^n}_{y_k=1} \chi_y\bigl( u(x^k)\bigr) +
\mathop{\sum_{k=j+1}^n}_{y_k=-1} \chi_y\bigl( u(x^k)\bigr),
$$
which is a sum of two geometric series, where the first contains
$m_j(y)$ elements and the second contains $n-j-m_j(y)$ elements.
Hence, this equals
$$
\sum_{l=0}^{m_j(y)-1} pq^l + \sum_{l=0}^{n-j-m_j(y)-1} -q^l  =
p{{1-q^{m_j(y)}}\over{1-q}} - {{1-q^{n-j-m_j(y)}}\over{1-q}}.
$$
and consequently
$$
p-1 + (q-1) \sum_{k=j+1}^n \chi_y\bigl(u(x^k)\bigr) =
pq^{m_j(y)} - q^{n-j-m_j(y)}.
$$
Thus the quadratic relation is preserved if
$\chi=\chi_y\bigl(u(x^j)\bigr)$ satisfies
$$
\chi^2 = \bigl( pq^{m_j(y)} - q^{n-j-m_j(y)}\bigr)\chi + pq^{n-j}
$$
or $\chi=pq^{m_j(y)}$ or $\chi=-q^{n-j-m_j(y)}$.
\qed\enddemo

\demo{Remark \theoremname{\remthmcharactersV}} Observe that
two of these characters can be easily calculated for arbitrary
$u(x)\in V_n$;
$$
\align
\chi_{(1,\ldots,1)}\bigl(u(x)\bigr) &=\iota(T_{u_x}) =
p^{w(x)} q^{nw(x) - \sum_{j=1}^{w(x)} i_j}, \\
\chi_{(-1,\ldots,-1)}\bigl(u(x) \bigr) &= \iota^\prime(T_{u_x})
=(-1)^{w(x)} q^{nw(x) - \sum_{j=1}^{w(x)} i_j}.
\endalign
$$
The characters $\chi_{(1,\ldots,1)}$, respectively
$\chi_{(-1,\ldots,-1)}$,
are the images of the one-dimensional representations $\iota$,
respectively $\iota^\prime$, of ${\Cal H}_n$ under the projection
$\p\colon {\Cal H}_n \to V_n$.
\enddemo

The algebra $V_n$ is a split semisimple algebra over
$\C(p^{1/2},q^{1/2})$, since it is a commutative $2^n$-dimensional
algebra with $2^n$ different one-dimensional representations. So we
can apply Kilmoyer's results \cite{\CurtRone, \S 9B} to obtain part
of the following orthogonality relations.

\proclaim{Proposition \theoremname{\proporthocharacters}}
The characters $\chi_y$, $y\in\Ztn$, of $V_n$ satisfy the
following orthogonality relations; for $y,z\in\Ztn$
$$
\sum_{x\in \Ztn} \chi_y\bigl( \hat u(x)\bigr) \,
\chi_z\bigl( \hat u(x)\bigr) = \d_{y,z} h_y, \quad
h_y = \prod_{j=1}^n 1+ \bigl(pq^{2m_j(y)+j-n}\bigr)^{y_j}.
$$
\endproclaim

\proclaim{Corollary \theoremname{\corproporthocharacters}}
The dual orthogonality relations hold; for $x,z\in\Ztn$
$$
\sum_{y\in \Ztn} {1\over{h_y}}
\chi_y\bigl( \hat u(x)\bigr) \, \chi_y\bigl( \hat u(z)\bigr) =
\d_{x,z}.
$$
\endproclaim

\demo{Proof of Proposition \thmref{\proporthocharacters}}
Apart from the squared norm the proposition follows immediately
from \cite{\CurtRone, Prop.~(9.17), (9.19)}.
To calculate the squared norm we observe that
$$
h_y = \sum_{x\in\Ztn} {{\Bigl(\chi_y\bigl(u(x)\bigr)\Bigr)^2}\over
{\iota(T_{u_x})}} =
\Bigl( 1 +
\Bigl(\chi_y\bigl(u(x^1)\bigr)\Bigr)^2 p^{-1}q^{1-n}\Bigr)
\sum_{z\in\Zt^{n-1}} {{\Bigl(\chi_y\bigl(u(1,z)\bigr)\Bigr)^2}\over
{\iota(T_{u(1,z)})}}
$$
by writing $x\in\Ztn$ as $(1,z)$ and $(-1,z)$ for $z\in\Zt^{n-1}$ and
using \thetag{\vglontbindingbasis}, $\chi_y$ being a character and
Proposition \thmref{\propdistcosetreprs}. Note that, by
\thetag{\vglontbindingbasis} and Theorem \thmref{\thmcharactersV},
the sum is independent of $y_1$. Theorem \thmref{\thmcharactersV}
implies that
$$
\Bigl(\chi_y\bigl(u(x^1)\bigr)\Bigr)^2 p^{-1}q^{1-n} =
\bigl( pq^{2m_1(y)+1-n}\bigr)^{y_1}.
$$
Now iterate this procedure to find the value for $h_y$. \qed
\enddemo

\demo{Remark \theoremname{\remproporthocharacters}}
By Remark \thmref{\remthmcharactersV} we have $h_{(1,\ldots,1)}=
\sum_{x\in\Ztn}\iota(T_{u_x})$. In this case we also have,
by Proposition \thmref{\propdistcosetreprs},
$$
h_{(1,\ldots,1)} P_A(q) = \sum_{x\in\Ztn} \iota(T_{u_x})
\sum_{\s\in S_n} \iota(T_\s) = \sum_{w\in H_n} \iota(T_w) = P_B(p,q),
$$
where $P_B(p,q)$ is the Poincar\'e polynomial for the hyperoctahedral
group, which is defined by the last equality. Hence,
$h_{(1,\ldots,1)}=P_B(p,q)/P_A(q) = (-p;q)_n$ which can be checked
directly from the explicit expressions for the Poincar\'e
polynomials;
$$
P_A(q) =  {{(q;q)_n}\over{(1-q)^n}},
\qquad P_B(p,q) = {{ (-p;q)_n(q;q)_n}\over{(1-q)^n}},
$$
see \cite{\Macd}.

Similarly, since $\iota(T_\s)=\iota^\prime(T_\s)$, $\s\in S_n$, and
$\iota(T_n)=p$, $\iota^\prime(T_n)=-1$, we get
$$
h_{(-1,\ldots,-1)} P_A(q) = \sum_{x\in\Ztn}
{{\bigl(\iota^\prime(T_{u_x}) \bigr)^2}
\over{\iota(T_{u_x})}}
\sum_{\s\in S_n} {{\bigl(\iota^\prime(T_\s) \bigr)^2}
\over{\iota(T_\s)}} = \sum_{w\in H_n}
{{\bigl(\iota^\prime(T_w) \bigr)^2}
\over{\iota(T_w)}}
= P_B(p^{-1},q).
$$
Hence, $h_{(-1,\ldots,-1)}=P_B(p^{-1},q)/P_A(q)=(-p^{-1};q)_n$.
\enddemo

\demo{Remark \theoremname{\remprimitiveidemp}} Since we have a
non-degenerate bilinear form $B$ on $V_n$, see Theorem
\thmref{\propalgebrastructureV}(iii), we can identify
$\chi_y\in V_n^\ast$ with $\xi_y\in V_n$ by
$$
\chi_y(v) = B(v,\xi_y), \qquad \forall\, v\in V_n.
$$
Then Proposition \thmref{\proporthocharacters} states that
$\{\xi_y\mid y\in\Ztn\}$ is an orthogonal basis of $V_n$, or
$B(\xi_y,\xi_z)=\d_{y,z} h_y$. From \cite{\CurtRone, Prop.~(9.17)}
it follows that $h_y^{-1}\xi_y$ is the (central) primitive idempotent
corresponding to the one-dimensional representation $\chi_y$ of $V_n$.
Hence, we have
$$
\xi_y\xi_z = \d_{y,z} h_y \xi_y, \qquad \sum_{y\in\Ztn}
h_y^{-1} \xi_y = u(1,\ldots,1)
$$
as identities in $V_n$, since $u(1,\ldots,1)$ is the identity
of $V_n$.

More general, let the $\C(p^{1/2},q^{1/2})$-linear isomorphism
$b\colon V_n^\ast\to V_n$, $b(f) =\xi_f$, be defined
by $f(v)=B(v,\xi_f)$ for all $v\in V_n$. Then $b$ is the analogue
of the Fourier transform. To see this define
$\hat{\ }\colon V_n^\ast \to V_n^{\ast\ast}$ by
$$
\hat f (\chi) = \sum_{x\in\Ztn} f\bigl(\hat u(x)\bigr)
\chi\bigl(\hat u(x)\bigr), \qquad \chi\in V_n^\ast,
$$
and let $a\colon V_n^{\ast\ast}\to V_n$ be the standard
isomorphism, i.e. for $\eta\in V_n^{\ast\ast}$ we have
$\eta(\chi)=\chi\bigl(a(\eta)\bigr)$ for all $\chi\in V_n^{\ast}$.
Then $a(\hat f)=\xi_f=b(f)$.

Using the Fourier transform $b$ we can put a commutative
$\C(p^{1/2},q^{1/2})$-algebra structure on $V_n^\ast$
by $\chi\star\eta = b^{-1}\bigl( b(\chi)b(\eta)\bigr)$, which is an
analogue of the convolution product.
Since $b(\chi_y)$,
$y\in\Ztn$, are primitive idempotents we get the following analogue
of the convolution product of characters;
$$
\chi_y \star \chi_z = \d_{y,z} h_y\, \chi_y, \qquad y,z\in\Ztn.
$$
Also $\t = b^{-1}\bigl(u(1,\ldots,1)\bigr) =
\sum_{y\in\Ztn} h_y^{-1}\chi_y$ is the identity of the algebra
$V_n^\ast$ given in Theorem \thmref{\propalgebrastructureV}(iii).
\enddemo


\subhead\newsection Contragredient representations
\endsubhead

We can define the contragredient representation $\r^\ast$
of ${\Cal H}_n$ in $V_n^\ast$ as follows; for $\chi\in V_n^\ast$ put
$$
\bigl( \r^\ast (\sum_{w\in H_n}c_w\,T_w) \chi_y\bigr)(v) =
\chi_y\bigl( \r((\sum_{w\in H_n} c_w\, T_w)^\ast) v\bigr),
\qquad v\in V_n,
$$
for some antimultiplicative map
$\ast\colon{\Cal H}_n\to{\Cal H}_n$
which also preserves the quadratic relations of
\thetag{\vgldefHeckeBn}. There are two obvious candidates for
the $\ast$-operator;
$\ast_1$ is the antilinear map defined by
$T_w^{\ast_1}=T_{w^{-1}}$,
$(f(p^{1/2},q^{1/2}))^{\ast_1}=\overline{f(p^{1/2},q^{1/2})}$
for $f\in\C(p^{1/2},q^{1/2})$, 
and $\ast_2$ is the antilinear map defined
by $T_w^{\ast_2}=T_w^{-1}$ and
$(f(p^{1/2},q^{1/2}))^{\ast_2}=\overline{f(p^{-1/2},q^{-1/2})}$
for $f\in\C(p^{1/2},q^{1/2})$. Observe that these maps are
involutions and that $\ast_1$ and $\ast_2$ commute. The composition
$\ast_1\circ\ast_2\colon{\Cal H}_n\to{\Cal H}_n$ is the
multiplicative involution which gives rise to the Kazhdan-Lusztig
basis of ${\Cal H}_n$, cf. 
\cite{\Hump, Ch.~7}.
Let $\r^\ast$ be the contragredient
representation defined using $\ast_1=\ast$;
the contragredient representation $\r_0^\ast$
using $\ast_2$ is related to it, cf.
Remark 6.2. 

Since the basis $\hat u(x)$, $x\in\Ztn$, is orthonormal with
respect to the bilinear form $B$ of Theorem
\thmref{\propalgebrastructureV}(iii)
and $\r(T_i)$, $1\leq i<n$, is given by the orthogonal matrix $R$ in
\thetag{\vglRmatrixTi},
we see that $B(\r(T_w)v,w)=B(v,\r(T_w^\ast)w)$
for $w\in S_n$.
Since $\rho(T_n)$ is also given by an orthogonal matrix with
respect to the orthonormal basis $\{ \hat u(x)\}$ we see that
this holds for arbitrary $w\in H_n$.
This implies the following commutation relations for the action
of ${\Cal H}_n$ with the Fourier transform $b$ introduced in
Remark \thmref{\remprimitiveidemp};
$$
b\bigl( \r^\ast(T) \chi) = \r(T)\, b(\chi), \qquad
b\bigl( \r_0^\ast(T) \chi) = \r((T^{\ast_1})^{\ast_2})\, b(\chi),
\tag\eqname{\vglactierhosterconvo}
$$
for $\chi\in V_n^\ast$, $T\in{\Cal H}_n$.

\proclaim{Theorem \theoremname{\thmcontragrepone}} The action
of the generators of ${\Cal H}_n$ on the characters
$\chi_y\in V_n^\ast$ in $\r^\ast$ is given by
$$
\r^\ast(T_n)\, \chi_y = \chi_y \cases
p, &\text{if $y_n=1$,} \\
-1, &\text{if $y_n=-1$,}\endcases
$$
and for $1\leq i<n$ by
$$
\align
\r^\ast(T_i)\, \chi_y &=
q\chi_y, \qquad \text{if}\  y_i=y_{i+1}, \\
\r^\ast(T_i)\, \chi_y &=
{{pq^{m_i(y)}(q-1) \chi_y +
(pq^{m_i(y)}+q^{n-i-m_i(y)})\chi_{y^{s_i}} }
\over{pq^{m_i(y)}+q^{n-i-1-m_i(y)}}},
\qquad \text{if}\  y_i=1,\, y_{i+1}=-1, \\
\r^\ast(T_i)\, \chi_y &=
{{q^{n-i-m_i(y)}(q-1) \chi_y +
(pq^{m_i(y)}+q^{n-i-m_i(y)})\chi_{y^{s_i}} }
\over{pq^{m_i(y)-1}+q^{n-i-m_i(y)}}},
\qquad \text{if}\  y_i=-1,\, y_{i+1}=1,
\endalign
$$
where $y^\s$ is defined by \thetag{\vglcommrelZnSn} and $m_i(y)$
is as in Theorem \thmref{\thmcharactersV}.
\endproclaim

\demo{Remark \theoremname{\remthmcontragrepone}} Note that
for $1\leq i<n$ the result of Theorem \thmref{\thmcontragrepone}
can be written uniformly as
$$
\r^\ast(T_i)\chi_y =
{{y_i \chi_y\bigl((\p(T_{u(x^i)})\bigr)(q-1) \chi_y +
(pq^{m_i(y)}+q^{n-i-m_i(y)})\chi_{y^{s_i}} }
\over{pq^{m_{i+1}(y)}+q^{n-i-1-m_{i+1}(y)}}}.
$$
If we write this as $\r^\ast(T_i)\chi_y =A_i(y)\chi_y +
B_i(y)\chi_{y^{s_i}}$,
then $\r^\ast(T_i-q)\chi_y = B_i(y)(\chi_{y^{s_i}}-\chi_y)$.
The action of $\r_0^\ast(T_i)$ on $\chi_y$ also follows, since by
\thetag{\vglactierhosterconvo}
$$
\r_0^\ast(T_i)\chi_y = b^{-1}\bigl(\r((T_i^{\ast_1})^{\ast_2})
b(\chi_y)\bigr) =b^{-1}\bigl(\r(T_i^{-1}) b(\chi_y)\bigr) =
\r^\ast(T_i^{-1})\chi_y.
$$
Now use $T_i^{-1}=q^{-1}T_i+q^{-1}-1$ ($1\leq i<n$) and a
similar expression for $T_n^{-1}$ with $q$ replaced by $p$.
\enddemo

\demo{Proof} Let us first consider the action of $\r^\ast(T_n)$,
then Theorem \thmref{\propalgebrastructureV},
\thetag{\vglontbindingbasis} and $\chi_y$ being a character
we obtain
$$
\bigl(\r^\ast(T_n)\chi_y\bigr) \bigl( u(x)\bigr)  =
\chi_y\bigl( u(x)\bigr) \cases
\chi_y\bigl( u(x^n)\bigr), &\text{if $x_n=1$,}\\
(p-1) + p \bigl( \chi_y\bigl( u(x^n)\bigr)\bigr)^{-1},
&\text{if $x_n=-1$.} \endcases
$$
Now by Theorem \thmref{\thmcharactersV}
$\chi_y\bigl( u(x^n)\bigr)$ equals $-1$ if $y_n=-1$ and $p$ if
$y_n=1$, which implies the result for the action of $T_n$.

Theorem \thmref{\propalgebrastructureV} gives for $0\leq i<n$
$$
\bigl(\r^\ast(T_i)\chi_y\bigr) \bigl( u(x)\bigr)  =
\cases \chi_y\bigl(u(x^{s_i})\bigr), &\text{if $x_i=1$,
$x_{i+1}=-1$,}\\
q\chi_y\bigl(u(x)\bigr), &\text{if $x_i=x_{i+1}$,}\\
(q-1)\chi_y\bigl(u(x)\bigr) +q \chi_y\bigl(u(x^{s_i})\bigr),
&\text{if $x_i=-1$, $x_{i+1}=1$.}\endcases
\tag\eqname{\vglhulpcontrag}
$$
In case $y_i=y_{i+1}$ we have $\chi_y\bigl(u(x^{s_i})\bigr)=
q\chi_y\bigl(u(x)\bigr)$ if $x_i=1$, $x_{i+1}=-1$. This follows
from \thetag{\vglontbindingbasis}
and $\chi_y$ being a character and
$\chi_y\bigl(u(x^i)\bigr)=
q\chi_y\bigl(u(x^{i+1})\bigr)$, which follows from Theorem
\thmref{\thmcharactersV}. Using this in \thetag{\vglhulpcontrag}
shows $\r^\ast(T_i)\chi_y = q\chi_y$ for $y_i=y_{i+1}$.

We now consider the case $y_i=1$, $y_{i+1}=-1$, so we have to show
that \thetag{\vglhulpcontrag} equals
$$
{{pq^{m_i(y)}(q-1) \chi_y\bigl( u(x)\bigr)
 + (pq^{m_i(y)}+q^{n-i-m_i(y)})\chi_{y^{s_i}}\bigl( u(x)\bigr) }
\over{pq^{m_i(y)}+q^{n-i-1-m_i(y)}}}
\tag\eqname{\vglhulptweecontrag}
$$
for all choices of $x$. To treat the case $x_i=x_{i+1}$
we note that $\chi_y\bigl( u(x)\bigr) = \chi_{y^{s_i}}
\bigl( u(x)\bigr)$, which follows from Theorem
\thmref{\thmcharactersV} and $m_p(y^{s_i})=m_p(y)$ for $p\not= i$
and $m_i(y^{s_i})=m_i(y)+1$. This shows that for $x_i=x_{i+1}$
\thetag{\vglhulptweecontrag} reduces to $q\chi_y\bigl(u(x)\bigr)$.

To treat the case $x_i\not= x_{i+1}$ we introduce $z\in\Ztn$ defined
by $z_p=x_p$ for $p\not= i,i+1$ and $z_i=z_{i+1}=1$. From
the previous paragraph we get
$\chi_y\bigl(u(z)\bigr) = \chi_{y^{s_i}}
\bigl( u(z)\bigr)$. Using Theorem
\thmref{\thmcharactersV} and $\chi_y$ being a character shows that
\thetag{\vglhulpcontrag} equals
$$
\chi_y \bigl( u(z)\bigr)
\cases pq^{m_i(y)}, &\text{if $x_i=1$, $x_{i+1}=-1$,}\\
(q-1)pq^{m_i(y)} - q^{n-i-m_i(y)},
&\text{if $x_i=-1$, $x_{i+1}=1$,}\endcases
\tag\eqname{\vglhuldriecontrag}
$$
since $m_i(y)=m_{i+1}(y)$ in this case. Because, by
Theorem \thmref{\thmcharactersV},
$$
\chi_y \bigl( u(x)\bigr) =
\chi_y \bigl( u(z)\bigr) \cases
-q^{n-i-1-m_i(y)}, &\text{if $x_i=1$, $x_{i+1}=-1$,}\\
pq^{m_i(y)}&\text{if $x_i=-1$, $x_{i+1}=1$,}\endcases
$$
and
$$
\chi_{y^{s_i}} \bigl( u(x)\bigr) =
\chi_y \bigl( u(z)\bigr) \cases
pq^{m_i(y)}, &\text{if $x_i=1$, $x_{i+1}=-1$,}\\
-q^{n-i-m_i(y)-1}&\text{if $x_i=-1$, $x_{i+1}=1$,}\endcases
$$
it is straightforward to check that \thetag{\vglhulptweecontrag}
equals \thetag{\vglhuldriecontrag}. This gives the
action of $T_i$ on $\chi_y$ with $y_i=1$, $y_{i+1}=-1$.

To prove the last identity of the theorem, we apply $\r^\ast(T_i)$
once more to $\r^\ast(T_i)\chi_y$ with $y_i=1$, $y_{i+1}=-1$.
Using the quadratic relation for $T_i$ we get an expression
for $\r^\ast(T_i)\chi_{y^{s_i}}$ in terms of $\chi_y$
and $\r^\ast(T_i)\chi_y$. Using the result already proved
for the last term gives an expression for
$\r^\ast(T_i)\chi_{y^{s_i}}$ in terms of $\chi_y$
and $\chi_{y^{s_i}}$. Now replace $y$ by $y^{s_i}$ to find the
result after some calculations. \qed
\enddemo

\proclaim{Corollary \theoremname{\corthmcontragrepone}}
Define for $f=0,\ldots,n$ the space $U_f$ of
characters $\chi_y$ with $w(y)=f$;
$$
U_f = \bigoplus_{w(y)=f} \C(p^{1/2},q^{1/2})\,
\chi_y \subset V_n^\ast,
$$
so that $\dim_{\C(p^{1/2},q^{1/2})} U_f = {n\choose f}$.
Then $V_n^\ast = \bigoplus_{f=0}^n U_f$
is the decomposition of the representation $\r^\ast$ in $V_n^\ast$
into irreducible ${\Cal H}_n$-modules.
\endproclaim

\demo{Remark \theoremname{\remlinkHoefsmit}} The irreducible
representations of ${\Cal H}_n$ have been classified by Hoefsmit
\cite{\Hoef, Def. (2.2.6), Thm. (2.2.7), (2.2.14)}
and are parametrised by double partitions of $n$. The irreducible
${\Cal H}_n$-module $U_f$ corresponds to $\bigl( (n-f), (f)\bigr)$
and $\chi_y$ corresponds to the standard tableau of shape
$\bigl( (n-f), (f)\bigr)$ given by
$$
\alignat2
&\boxed{n+1-j_{n-f}\mid n+1-j_{n-f-1} \mid \ldots
\mid n+1-j_1} &&\quad (n-f), \\
&\boxed{n+1-i_f\mid n+1-i_{f-1} \mid \ldots \mid n+1-i_1}
&&\quad (f),
\endalignat
$$
where $y_{i_1}=\ldots=y_{i_f}=-1$, $i_1<\ldots<i_f$, and
$y_{j_1}=\ldots=y_{j_{n-f}}=1$, $j_1<\ldots<j_{n-f}$.
Applying \cite{\Hoef, Prop.~(3.3.3)} shows that for $x\in\Ztn$
the operator $\r^\ast(T_x)$ is diagonal with respect to the
basis $\{ \chi_y\}_{y\in\Ztn}$ of $V_n^\ast$.
\enddemo

As before we call $\chi\in V_n^\ast$ a ${\Cal F}_n$-invariant element
of $V_n^\ast$ if $\r^\ast(T_\s)\, \chi=q^{\ell(\s)}\, \chi$
for all $\s\in S_n$. So $\chi$ is ${\Cal F}_n$-invariant if it
realises the index representation $\iota$ of ${\Cal F}_n$ in
$V_n^\ast$.

Let $V_n^d = \bigoplus_{\{ x\in\Ztn\mid w(x)=d\} }
\C(p^{1/2},q^{1/2})\, v(x) \subset V_n$,
then $V_n =\bigoplus_{d=0}^n V_n^d$ and by \thetag{\vglTiopvx} it
follows that each $V_n^d$ is invariant under the action of
${\Cal F}_n$ via $\r$. In general $V_n^d$ is not an irreducible
${\Cal F}_n$-module, and it can be obtained from inducing the
index representation of ${\Cal F}_d\otimes {\Cal F}_{n-d}$. For
the decomposition as ${\Cal F}_n$-module we can proceed
as in Dunkl \cite{\Dunk, \S 2}.

We can now describe the $\U$-module of ${\Cal F}_n$-invariant
elements in $V_n$ explicitly.

\proclaim{Theorem \theoremname{\thmexplictcdescFinveltsinV}}
Define the ${\Cal F}_n$-invariant elements $w_d$, $d=0,1,\ldots,n$,
in $V_n$ by
$$
w_d =\r(P)\,
v({\undersetbrace n-d \to{1,\ldots,1}},
{\undersetbrace d \to{-1,\ldots,-1}}) \in (V_n)^{{\Cal F}_n},
$$
with $P$ and $v(x)$ defined in \thetag{\vgldefP} and
Theorem \thmref{\propalgebrastructureV}(iii),
then $w_d$ is non-zero and forms an orthogonal
basis for the space of ${\Cal F}_n$-invariant
elements in $V_n$;
$$
B(w_d,w_e) = \d_{d,e}\, p^{-d}q^{-d(d-1)/2}
\left[ {n\atop d}\right]_q^{-1}.
$$
Moreover, for $x\in\Ztn$ with $w(x)=d$ we have
$\r(P)\, v(x) = w_d$.
The action of $\U$ in this basis is given by
$t(K)\, w_d = q^{d-n/2}\, w_d$ and
$$
t(E)\, w_d = p^{1/2}q^{(1-n)/2} q^d {{1-q^{n-d}}\over{1-q}}\,
w_{d+1}, \quad
t(F)\, w_d = p^{-1/2} q^{1-d} {{1-q^d}\over{1-q}}\, w_{d-1}.
$$
\endproclaim

Here $\dsize{\left[ {n\atop d}\right]_q=
{{(q^n;q^{-1})_d}\over{(q;q)_d}}}$ denotes the $q$-binomial coefficient 

\demo{Proof} From \thetag{\vglDeltanopAB},
\thetag{\vglactiefundrep} and
the identification \thetag{\vglidenWnenVn} we see that
$$
t(K) \colon V_n^d \to V_n^d, \quad
t(E) \colon V_n^d \to V_n^{d+1}, \quad
t(F) \colon V_n^d \to V_n^{d-1},
$$
using a similar expression for $\Delta^{(n)}(F)$. From Theorem
\thmref{\thmreprquantumsl2} and
Proposition \thmref{\propFinvariantelts}
it follows that the space of ${\Cal F}_n$-invariant elements in
$V_n^d$ is one-dimensional.

Theorem \thmref{\thmcontragrepone} implies that
$\chi_1=\chi_{(1,\ldots,1)}$
is ${\Cal F}_n$-invariant in $V_n^\ast$ under $\r^\ast$. Using
Remark \thmref{\remthmcharactersV} we get
$$
1 = \chi_1\bigl( v(x)\bigr) =
\Bigl( \r^\ast(P)\chi_1\Bigl) \bigl(v(x)\bigr)
= \chi_1 \Bigl( \r(P)\, v(x)\Bigr).
$$
Since the left hand side is non-zero it follows that Hecke
symmetrising any basis elements yields a non-zero
${\Cal F}_n$-invariant element
in $V_n^{w(x)}$. Thus $\r(P)\, v(x)
=C(x)\, w_{w(x)}$ for some non-zero $C(x)\in\C(p^{1/2},q^{1/2})$.
So we obtain $1= C(x)\, \chi_1(w_{w(x)})$.
Since by definition $C(1,\ldots,1,-1,\ldots,-1)=1$ we have
$\chi_1(w_d)= 1$ and $C(x) = 1/\chi_1(w_{w(x)}) = 1$.

To calculate the action of $K$ we note that by
\thetag{\vglactiefundrep} and \thetag{\vglDeltanopAB}
$$
t(K)\, \hat u (1,\ldots,1,-1,\ldots,-1) =
q^{-(n-d)/2} q^{d/2}\, \hat u (1,\ldots,1,-1,\ldots,-1),
$$
in terms of  the orthonormal basis $\hat u(x)$ of $V_n$.
Since $t(K)$ commutes with the Hecke symmetrisator by Jimbo's
Theorem \thmref{\thmJimbo} we find $t(K)w_d=q^{d-n/2}w_d$.

Similarly,
$$
t(E)\, \hat u (1,\ldots,1,-1,\ldots,-1) =
\sum_{i=1}^{n-d} q^{-(i-1)/2} \,
\hat u (1,\ldots,1,-1,1,\ldots,1,-1,\ldots,-1)
$$
where the first $-1$ occurs at the $i$-th place. Now apply the
Hecke symmetrisator, use Jimbo's Theorem \thmref{\thmJimbo},
and $v(x)=\iota(T_{u_x})^{-1/2}\hat u(x)$ to find
$$
t(E)\, w_d = \Bigl(\sum_{i=1}^{n-d} q^{-(i-1)/2}
p^{1/2} q^{(n-i)/2} \Bigr)\, w_{d+1},
$$
which gives the result.

The action of $t(F)$ can be calculated similarly, but can
also be derived from the commutation relation for $E$ and $F$ in
$\U$, see Definition \thmref{\defquantumsl2}. So if
$t(F)w_d=c_dw_{d-1}$ we obtain
$$
p^{1/2}q^{d-1}{{q^{(1-n)/2}}\over{1-q}} \bigl(
c_d (1-q^{n-d+1}) - c_{d+1} q (1-q^{n-d})\bigr)
= {{q^{d-n/2}-q^{-d+n/2}}\over{q^{1/2}-q^{-1/2}}}.
$$
Together with the initial condition $c_0=0$ this two-term recurrence
relation determines $c_d$.

Using the fact that $V_n$ is a $\ast$-module of $\U$ we get
$$
\multline
p^{1/2}q^{(1-n)/2} q^d {{1-q^{n-d}}\over{1-q}}\, B(w_{d+1},w_{d+1})
= B(t(E)w_d,w_{d+1}) = \\
B(w_d,t(q^{-1/2}FK)w_{d+1}) =
p^{-1/2} q^{(1-n)/2} {{1-q^{d+1}}\over{1-q}} B(w_d,w_d).
\endmultline
$$
Solve this recurrence relation with the condition $B(w_0,w_0)=1$
to find the result.
\qed\enddemo

\proclaim{Proposition \theoremname{\propJnboinVnstar}}
The space of ${\Cal F}_n$-invariant elements in $V_n^\ast$ is
$n+1$-dimensional, and each $U_f$ contains a one-dimensional
space of ${\Cal F}_n$-invariant elements. Moreover,
$\phi_f=\r^\ast(P)\, \chi_y = \r^\ast(P)\, \chi_z$
is a non-zero ${\Cal F}_n$-invariant element in $U_f$
for $y,z\in\Ztn$ with $w(y)=w(z)=f$.
\endproclaim

\demo{Remark \theoremname{\rempropJnboinVnstar}} Extend
$\phi_f\colon V_n\to \C(p^{1/2},q^{1/2})$ to ${\Cal H}_n$ by
defining $\phi_f(T)=\phi\bigl(\p(T)\bigr)$ for $T\in{\Cal H}_n$
and $\p\colon {\Cal H}_n\to V_n$ as in \S 3. 
Then $\phi_f$ is
left and right $\iota|_{{\Cal F}_n}$-invariant;
$\phi_f(T_\s T T_\t)= \iota(T_\s)\iota(T_\t) \phi_f(T)$
for all $\s,\t\in S_n$. Moreover, $\phi_f$ is contained in the
irreducible ${\Cal H}_n$-module $U_f$, so that we may regard
$\phi_f$ as a `zonal spherical function' on the Hecke algebra
${\Cal H}_n$.
\enddemo

\demo{Proof} First observe that if $\chi$ is
${\Cal F}_n$-invariant, then
$$
\chi\bigl( v(x)\bigr) = \Bigl(
\r^\ast(P)\, \chi\Bigr) \bigl(v(x)\bigr) =
\chi \Bigl(
\r(P)\, v(x) \Bigr) = \chi(w_{w(x)}),
$$
by Theorem \thmref{\thmexplictcdescFinveltsinV}.
So $\chi$ is completely determined by $\chi(w_d)$, $d=0,1,\ldots,n$.
The dimension of the space of ${\Cal F}_n$-invariant elements in
$V_n^\ast$ is at most $n+1$. If we can show that $U_d$ has at
least a one-dimensional subspace of ${\Cal F}_n$-invariant
elements, it follows that $U_d$ has a subspace of
${\Cal F}_n$-invariant elements of dimension one and the
space of ${\Cal F}_n$-invariant elements in $V_n^\ast$ is
$n+1$-dimensional.

Obviously, $\r^\ast(P)\, \chi_y$
is a ${\Cal F}_n$-invariant element in $U_{d(y)}$. Now
$$
\r^\ast(P)\, \chi_y \bigl( u(1,\ldots,1)\bigr) =
\chi_y\bigl( u(1,\ldots,1)\bigr) = 1,
$$
since the unit element $u(1,\ldots,1)$ of $V_n$ is
${\Cal F}_n$-invariant. So this is a non-zero
${\Cal F}_n$-invariant element in $U_f$, $f=w(y)$,
and by the previous paragraph the same element with $y$ replaced
by $z$ with $w(z)=f$ is a multiple of this
${\Cal F}_n$-invariant element in $U_f$.
The constant is $1$, since evaluated at
$u(1,\ldots,1)$ gives $1$ in both cases.
\qed \enddemo

\demo{Remark \theoremname{\remdeftster}}
Since the representation $t$ of $\U$ acts in $V_n$, we can
also define contragredient representations of $\U$ in $V_n^\ast$
using antimultiplicative mappings of $\U$
to itself. An obvious candidate is the antipode of the
Hopf-algebra $\U$.
Another candidate is the $\ast$-operator defined by
\thetag{\vglatstU} We define a representation $t^\ast$
$$
\bigl( t^\ast(X) \chi\bigr)(v) =
\chi \bigl(t(X^\ast)\, v\bigr)
$$
for $X\in\U$, $\chi\in V_n^\ast$ and $v\in V_n$, and we do not
use the contragredient representation defined using the antipode.
Observe that $B(t(X)v,w)=B(v,t(X^\ast)w)$ implies
$b\bigl( t^\ast(X) \chi\bigr) = t(X)b(\chi)$, with $b$ the Fourier
transform defined in Remark \thmref{\remprimitiveidemp}.
\enddemo


\subhead\newsection  ${\Cal F}_n$-invariant elements and
$q$-Krawtchouk polynomials
\endsubhead

It follows from Proposition \thmref{\propJnboinVnstar} that
$\phi_f = \r^\ast(P)\chi_y\in \, (V_n^\ast)^{{\Cal F}_n}$
only depends on $f=w(y)$.
Then $\phi_f\bigl(u(x)\bigr)$ is determined by $d=w(x)$ by
Theorem \thmref{\thmexplictcdescFinveltsinV} and
$$
\phi_f(w_d) = \chi_y(w_d), \qquad w(y)=f.
$$
Proposition \thmref{\proporthocharacters} implies that the
$\phi_f(w_d)$ satisfy certain orthogonality relations for $d$
running through the finite set $\{0,1,\ldots,n\}$.
This is not sufficient
to determine $\phi_f(w_d)$
since we do not a priori know if $\phi_f(w_d)$ is
a polynomial, but we come back to the orthogonality
properties later.  Using the fact that $\phi_f$ is a sum of
characters we give in Theorem 7.1 
an explicit expression for a weighted average
of $\phi_f$ evaluated in $w_{d-1}$, $w_d$ and $w_{d+1}$,
which corresponds to the finite Laplacian in Stanton
\cite{\StanAJM, Def.~5.2}.

\proclaim{Theorem \theoremname{\thmLaplacevgl}} The values
$\phi_f(w_d)$ satisfy the recurrence relation for $0\leq f\leq n$
$$
\multline
\bigl( p(1-q^{n-f}) -(1-q^f)\bigr)\,\phi_f(w_d)  =
pq^d(1-q^{n-d})\,\phi_f(w_{d+1}) \\ + (1-q^d)(p-1)\,
\phi_f(w_d) + (1-q^d)\, \phi_f(w_{d-1})
\endmultline
$$
with initial conditions
$$
\phi_f(w_0) =1,\qquad \phi_f(w_n) = (-p)^{f-n}q^{f(f-n)}.
$$
\endproclaim

\proclaim{Corollary \theoremname{\corthmLaplacevgl}} The basis
$\phi_f$ of the $n+1$-dimensional irreducible $\U$-module
$(V_n^\ast)^{{\Cal F}_n}$ are eigenvectors of the following
self-adjoint operator;
$$
t^\ast\bigl(E+E^\ast + {{p^{1/2}-p^{-1/2}}\over{q^{1/2}-q^{-1/2}}}
(K-1)\bigr)\, \phi_f
= {{p^{1/2}q^{n/2-f}-p^{-1/2}q^{f-n/2} +p^{-1/2}-p^{1/2}}
\over{q^{1/2}-q^{-1/2}}} \, \phi_f.
$$
\endproclaim

\demo{Proof} Combine the results of Theorems \thmref{\thmLaplacevgl}
and \thmref{\thmexplictcdescFinveltsinV}
with the definition of $t^\ast$. \qed
\enddemo

\demo{Proof of Theorem \thmref{\thmLaplacevgl}}
The values $\phi_f(w_0)$ and $\phi_f(w_n)$ follow from Theorem
\thmref{\thmcharactersV} since $V_n^0$ and $V_n^n$ are
one-dimensional.

Let $\phi_f = \sum_{w(y)=f} c_f(y) \chi_y$, then for $v\in V_n$
and $w\in V_n^{{\Cal F}_n}$
$$
\phi_f(wv) = \sum_{w(y)=f} c_f(y) \chi_y(w) \chi_y(v) =
\chi_y(w) \phi_f(v) = \phi_f(w) \phi_f(v),
\tag\eqname{\vglabstractprodform}
$$
since $\chi_y$ is a character of $V_n$ and $\chi_y(w)$ only depends
on $w(y)=f$.

We now use this with $w= \sum_{i=1}^n u(x^i)$, which is
${\Cal F}_n$-invariant as easily follows from Theorem
\thmref{\propalgebrastructureV}(i), and
by Theorem \thmref{\thmexplictcdescFinveltsinV} this is a constant
multiple of $w_1$. In order to evaluate $\phi_f(w)$
we may take $y=(-1,\ldots,-1,1,\ldots,1)$ with $w(y)=f$ to get
$$
\phi_f(w) = \sum_{i=1}^n \chi_y\bigl( u(x^i)\bigr) =
\sum_{i=1}^f -q^{f-i} + \sum_{i=f+1}^n pq^{n-i}
= - {{1-q^f}\over{1-q}} + p {{1-q^{n-f}}\over{1-q}}.
$$
We take $v=v(1,\ldots,1,-1,\ldots,-1)$ with $d$ minus signs, so that
$\phi_f(v)=\phi_f(w_d)$.

So it remains to calculate $\phi_f(wv)$, and for this we need
an explicit expression for $\r(P)(vw)$
in terms of the basis $w_e$, $0\leq e\leq n$, for
${\Cal F}_n$-invariant elements of Theorem
\thmref{\thmexplictcdescFinveltsinV}. Now observe that
$$
\multline
p^d q^{d(d-1)/2}vw = u(1,\ldots,1,-1,\ldots,-1)\, w =
\sum_{i=1}^{n-d}  u(1,\ldots,1,{\overset i\to{-1}},1,\ldots,1,
{\undersetbrace d\to{-1,\ldots,-1}}) \\
+ \sum_{i=n-d+1}^n u(1,\ldots,1,
{\undersetbrace i-n-d+1 \to{-1,\ldots,-1}}, {\overset i\to{1}},
{\undersetbrace n-i \to{-1,\ldots,-1}}) \, \bigl(u(x^i)\bigr)^2.
\endmultline
$$
Using Theorem \thmref{\thmexplictcdescFinveltsinV} we see that Hecke
symmetrising the first sum yields a multiple of $w_{d+1}$; explicitly
$$
\sum_{i=1}^{n-d} p^{d+1}q^{d(d-1)/2 +n-i} \, w_{d+1} =
p^{d+1} q^{d(d-1)/2} q^d {{1-q^{n-d}}\over{1-q}}\, w_{d+1}.
$$
In order to treat the second sum we need to the following Lemma.

\proclaim{Lemma \theoremname{\lemcalforqdiff}} For $z\in\Zt^{i-1}$
we have in $V_n$
$$
\multline
u(z, {\overset i\to{1}},
{\undersetbrace n-i \to{-1,\ldots,-1}}) \, \bigl(u(x^i)\bigr)^2
= q^{n-i} (p-1)\, u(z,-1,\ldots,-1) + \\
pq^{n-i-1}(q-1) \sum_{j=i+1}^n  u(z,-1,\ldots,-1,{\overset j\to 1},
-1,\ldots,-1) + pq^{n-i}\, u(z,1,-1,\ldots,-1).
\endmultline
$$
\endproclaim

Assuming Lemma \thmref{\lemcalforqdiff} we see that the second sum
can be written as
$$
\multline
(p-1)\sum_{i=n-d+1}^n q^{n-i} \,
u(1,\ldots,1,{\undersetbrace d\to{-1,\ldots,-1}}) + \\
(q-1) \sum_{i=n-d+1}^n pq^{n-i-1} \sum_{j=i+1}^n
u(1,\ldots,1,{\undersetbrace j-n+d-1\to{-1,\ldots,-1}},
{\overset j\to 1}, {\undersetbrace n-j\to{-1,\ldots,-1}}) +\\
\sum_{i=n-d+1}^n pq^{n-i} u(1,\ldots,1,{\undersetbrace i-n+d-1\to
{-1,\ldots,-1}},{\overset i \to 1}, {\undersetbrace n-i\to
{-1,\ldots,-1}})
\endmultline
$$
Now we can use Theorem \thmref{\thmexplictcdescFinveltsinV} to
Hecke symmetrise each of these sums. Hecke symmetrising the first
sum gives
$$
p^d q^{d(d-1)/2}(p-1)\sum_{i=n-d+1}^n q^{n-i} \, w_d =
p^d q^{d(d-1)/2}(p-1) {{1-q^d}\over{1-q}}\, w_d,
$$
and Hecke symmetrising the last sum gives
$$
\sum_{i=n-d+1}^n pq^{n-i} p^{d-1} q^{d(d-1)/2-(n-i)}\, w_{d-1} =
p^d q^{d(d-1)/2}(n-d)\, w_{d-1}
$$
and Hecke symmetrising the double sum gives
$$
\multline
(q-1) \sum_{i=n-d+1}^n pq^{n-i-1} \sum_{j=i+1}^n p^{d-1}
q^{d(d-1)/2-(n-j)} \, w_{d-1}
= \\ (q-1)p^d  q^{d(d-1)/2} \sum_{i=n-d+1}^n
{{1-q^{n-i}}\over{1-q}} \, w_{d-1} = p^d  q^{d(d-1)/2}
\bigl( {{1-q^d}\over{1-q}}-(n-d)\bigr)\, w_{d-1}.
\endmultline
$$

Collecting all terms shows
$$
(1-q)\r(P) \bigl(wv\bigr)
= pq^d(1-q^{n-d})\, w_{d+1} + (1-q^d)(p-1)\, w_d
+(1-q^d)\, w_{d-1}.
$$
Now apply $\phi_f$ to get
$$
(1-q) \phi_f(wv) = pq^d(1-q^{n-d})\, \phi_f(w_{d+1})
+ (1-q^d)(p-1)\, \phi_f(w_d)
+ (1-q^d)\, \phi_f(w_{d-1}).
$$
Combine this with the first paragraphs of the proof to find the
result. \qed
\enddemo

\demo{Proof of Lemma \thmref{\lemcalforqdiff}}
First recall from Theorem \thmref{\propalgebrastructureV}(ii) that
$$
u(x^i)^2 = (q-1)\sum_{l=i+1}^n  u(x^ix^l) +
(p-1) v(x^i) + pq^{n-i}.
$$
This in particular proves the case $i=n$, and for the general case
we obtain
$$
\multline
u(z, {\overset i\to{1}},
{\undersetbrace n-i \to{-1,\ldots,-1}}) \, \bigl(u(x^i)\bigr)^2 =
(p-1)
u(z,-1,\ldots,-1) +\\ pq^{n-i}\,
u(z,{\overset i \to 1},-1,\ldots,-1) +
(q-1)\sum_{l=i+1}^n u(z,-1,\ldots,-1,{\overset l\to 1},-1\ldots,-1)
\bigl( u(x^l)\bigr)^2.
\endmultline
$$
Proceeding by downward induction on $i$ we can rewrite each term
in the last sum. Collecting the coefficients of each basis element
$u(x)$ then proves the Lemma. \qed
\enddemo

\demo{Remark \theoremname{\remproductformula}} The second-order
difference equation for $\phi_f(w_d)$ with respect to $d$ of
Theorem \thmref{\thmLaplacevgl} is the consequence of the
homomorphism property \thetag{\vglabstractprodform} for
$w=w_1$ (up to a constant) and (sufficiently) arbitrary $v$. So
we have calculated the coefficients in the expansion, or product
formula,
$$
\phi_f(w_k)\, \phi_f(w_d) = \sum_{l=0}^n c_l(k,d)\, \phi_f(w_l)
$$
for $k=1$. The coefficients $c_l(k,d)$ are determined by
$$
\r(P)\bigl( w_k w_d) =\sum_{l=0}^n c_l(k,d)\, w_l,
$$
or with $w_d$ replaced by $v$ as in the proof of Theorem
\thmref{\thmLaplacevgl}. It follows from Theorem
\thmref{\thmexplictcdescFinveltsinV} that
$$
c_l(k,d)\, B(w_l,w_l) = B(\r(P)\bigl( w_k w_d), w_l) =
B(w_kw_d,w_l) = \t(w_kw_dw_l),
$$
which is symmetric in $k$, $l$ and $d$, since $V_n$
is a commutative algebra. In general these constants seem hard
to calculate, but for $q=1$ they have been calculated
explicitly, cf. Remark 7.10. 
\enddemo

\demo{Remark \theoremname{\remconvoprop}} From
\thetag{\vglabstractprodform} we can derive a convolution property
for $\phi_f$, cf. Remark \thmref{\remprimitiveidemp}. Note that
$b(\phi_f)\in(V_n)^{{\Cal F}_n}$ by \thetag{\vglactierhosterconvo},
so that for all $v\in V_n$,
$$
(\phi_f\star\phi_g)(v)= \phi_f\bigl( b(\phi_g)v)=
\phi_f\bigl(b(\phi_g)\bigr) \, \phi_f(v).
$$
Since $\chi_y\star\chi_z =\d_{y,z} h_y\chi_y$, cf. Remark
\thmref{\remprimitiveidemp}, we see that
$\phi_f\star\phi_g=0$ for $f\not= g$.
So $\phi_f\bigl(b(\phi_g)\bigr)=0$ for $f\not= g$ and
$$
\phi_f\star\phi_g = \d_{f,g}\, \phi_f\bigl(b(\phi_f)\bigr)
\, \phi_f.
$$
The convolution property of $\phi_f =
\sum_{y\in\Ztn\colon w(y)=f} c_f(y)\, \chi_y$ can be used to
calculate $c_f(y)$, since it implies
$h_y\,\bigl(c_f(y)\bigr)^2=
\phi_f\bigl(b(\phi_f)\bigr)\, c_f(y)$,
so that $c_f(y)=0$ or
$c_f(y)=\phi_f\bigl(b(\phi_f)\bigr)\, h_y^{-1}$.
The first case is excluded by the proof of the next Proposition,
where $\phi_f\bigl(b(\phi_f)$ is also explictly given.
\enddemo

\proclaim{Proposition \theoremname{\proporthorelKfd}} The
following orthogonality relations hold for $0\leq d,e,f\leq n$;
$$
\gather
\sum_{d=0}^n p^d q^{d(d-1)/2}\left[ {n\atop d}\right]_q
\, \phi_g(w_d)\, \phi_f(w_d) =
\d_{g,f}\, H_f, \\
\sum_{f=0}^n {1\over{H_f}}\, \phi_f(w_g)\, \phi_f(w_d) = \d_{d,g}\,
 p^d q^{d(d-1)/2} \left[ {n\atop d}\right]_q^{-1}
\endgather
$$
with
$$
H_f= \left[ {n\atop f}\right]_q^{-1}
{{(-pq^{-f};q)_{n+1}}\over{pq^{n-f}+q^f}} q^{f(f+1)/2} p^{-f}.
$$
\endproclaim

\demo{Proof} The explicit value for the weights follows from the
observation $\phi_f\bigl(b(\phi_g)\bigr)=0$ for $f\not= g$ of
Remark \thmref{\remconvoprop}. Since
$b(\phi_g)\in (V_n)^{{\Cal F}_n}$
we develop it in the orthogonal basis $w_d$, $0\leq d\leq n$;
$$
b(\phi_g) = \sum_{d=0}^n {{B\bigl(b(\phi_g),w_d\bigr)}
\over{B(w_d,w_d)}}\, w_d =
\sum_{d=0}^n {{\phi_g(w_d)}\over{B(w_d,w_d)}}\, w_d.
$$
Hence, for $f\not= g$ we have
$$
0=\phi_f\bigl( b(\phi_g)\bigr) =
\sum_{d=0}^n
{{\phi_g(w_d)\phi_f(w_d)}\over{B(w_d,w_d)}},
$$
so that the value for the weights follows from Theorem
\thmref{\thmexplictcdescFinveltsinV} and
$H_f= \phi_f\bigl( b(\phi_f)\bigr) =
B\bigl(b(\phi_f),b(\phi_f)\bigr)$.

In order to calculate squared norm $H_f$ we first observe that
$H_f$ is non-zero for every $f$ and that
$$
H_0 = h_{(1,\ldots,1)}=(-p;q)_n,\qquad
H_n = h_{(-1,\ldots,-1)}=(-p^{-1};q)_n,
$$
cf. Proposition \thmref{\proporthocharacters} and Remark
\thmref{\remproporthocharacters}. By Corollary
\thmref{\corproporthocharacters}
$$
{1\over{H_f}} = \sum_{y\in\Ztn\colon w(y)=f} {1\over{h_y}}
$$
with $h_y$ as in Proposition \thmref{\proporthocharacters}.
Initially this holds up to a scalar independent of $f$ and $n$,
and by the initial conditions this scalar equals $1$.

We now use the notation $h_y(n)=h_y$ to stress the $n$-dependence.
{}From Proposition \thmref{\proporthocharacters} we see that for
$y=(1,z)$, $z\in\Zt^{n-1}$, with $f=w(y)=d(z)$
$$
h_y(n) = (1+pq^{2m_1(y)+1-n})\, h_z(n-1) =
(1+pq^{n-1-2f})\, h_z(n-1)
$$
and for $y=(-1,z)$, $z\in\Zt^{n-1}$, with $f=w(y)=d(z)+1$
$$
h_y(n) = (1+p^{-1}q^{n-1-2m_1(y)})\, h_z(n-1) =
(1+p^{-1}q^{2f-n-1})\, h_z(n-1).
$$
Writing the sum over $y\in\Ztn$ as two sums of the form
$y=(1,z)$ for $z\in\Zt^{n-1}$ and $y=(-1,z)$ for
$z\in\Zt^{n-1}$ gives the recurrence relation, with
$H_f=H_f(n)$,
$$
{1\over{H_f(n)}} = {1\over{(1+pq^{n-1-2f})H_f(n-1)}}
+ {1\over{(1+p^{-1}q^{2f-n-1})H_{f-1}(n-1)}}.
$$
This recurrence relation together with the initial
conditions completely determines $H_f(n)$.

Put $H_f(n)^{-1} = (pq^{n-f}+q^f)\, D_f^n$ to find the
recurrence relation
$$
(1+pq^{n-2f})\, D_f^n = D_f^{n-1} + pq^{n-2f}\, D_{f-1}^{n-1},
$$
which, by substition of $D_f^n =p^fq^{-f(f+1)/2}
(-pq^{-f};q)_{n+1}^{-1} E^f_n$, can be rephrased as
$$
(1+pq^{n-2f})\, E_f^n = (1+pq^{n-f})\, E_f^{n-1} +
(q^{n-f} + pq^{n-2f})\, E_{f-1}^{n-1}.
$$
The initial condition is now $E_0^n=1=E_n^n$. Comparing this
recurrence relation with
$$
\left[ {n\atop f}\right]_q = q^f \left[ {{n-1}\atop f}\right]_q
+ \left[ {{n-1}\atop {f-1}}\right]_q
= \left[ {{n-1}\atop f}\right]_q
+ q^{n-f} \left[ {{n-1}\atop {f-1}}\right]_q
$$
gives $E_f^n = \left[{n\atop f}\right]_q$. \qed
\enddemo

Theorem \thmref{\thmLaplacevgl}
completely determines the values of $\phi_f(w_d)$ and this can be
expressed using $q$-Krawtchouk polynomials.
Define the $q$-Krawtchouk polynomials of degree $n$,
$0\leq n\leq N$, and of argument $q^{-x}$ by
$$
K_n(q^{-x};a,N;q)  = {}_3\vp_2 \left(
{{q^{-n},q^{-x},-q^{n-N}/a}\atop {q^{-N},\ 0}};q,q\right),
\tag\eqname{\vgldefqKrawtchouk}
$$
see e.g. \cite{\StanAJM, \S 3}, \cite{\GaspR, ex.~7.8(i)}.
The following Proposition can be found in
Stanton \cite{\StanAJM, Prop.~3.7}.

\proclaim{Proposition \theoremname{\proprecursionKrawtchouk}}
The $q$-Krawtchouk polynomials satisfy the second order
$q$-difference equation
$$
\multline
(q^n-aq^{N-n})\, K_n(q^{-x};a,N;q) = aq^x(1-q^{N-x}) \,
K_n(q^{-(x+1)};a,N;q) \\ + q^x(a -1)\, K_n(q^{-x};a,N;q)
+ (1-q^x)\, K_n(q^{-(x-1)};a,N;q)
\endmultline
$$
with initial conditions
$$
K_n(1;a,N;q)= 1, \qquad K_n(q^{-N};a,N;q) = (-a)^{-n} q^{n(n-N)}.
$$
\endproclaim

\proclaim{Theorem \theoremname{\thmsphericalasqKrawtchouk}}
We have
$\phi_f(w_d) = K_f(q^{-d};p,n;q)$
with the $q$-Krawtchouk polynomial $K_f$ defined by
\thetag{\vgldefqKrawtchouk}.
\endproclaim

\demo{Proof} Compare
Theorem \thmref{\thmLaplacevgl}
with Proposition \thmref{\proprecursionKrawtchouk}. \qed
\enddemo

Now that we have the zonal spherical function $\phi_f$ in terms
of explicit polynomials, we can interpret some of the
identities derived for $\phi_f$ as identities for
$q$-Krawtchouk polynomials. First of all,
the orthogonality relations of Proposition
\thmref{\proporthorelKfd} correspond to the
orthogonality relations for the $q$-Krawtchouk polynomials
and for the dual $q$-Krawtchouk polynomials, where the dual
$q$-Krawtchouk polynomials are defined by
$$
R_n(q^{-x}-q^{x-N}/a;a,N;q) = K_x(q^{-n};a,N;q).
$$
The second order $q$-difference equation for the $q$-Krawtchouk
polynomials of Theorem \thmref{\thmLaplacevgl}
is the three-term recurrence relation for
the dual $q$-Krawtchouk polynomials.

Secondly, Corollary
\thmref{\corthmLaplacevgl} corresponds to the fact that the
matrix elements of the
transition of the basis of eigenvectors for the action of $K$
to the basis of eigenvectors for the action of
$E+E^\ast + {{p^{1/2}-p^{-1/2}}\over{q^{1/2}-q^{-1/2}}}
(K-1)$ is given by $q$-Krawtchouk polynomials, see
Koornwinder \cite{\KoorZSE, Thm.~4.3}. Using this interpretation
Koornwinder \cite{\KoorZSE} is able to give an interpretation
of Askey-Wilson polynomials on the quantum $SU(2)$ group
as zonal spherical functions.

\demo{Remark \theoremname{\remalternativeproof}}
$V_{n-1}$ can be viewed as a subalgebra of $V_n$ by identifying
$u(x)\in V_{n-1}$ for $x\in\Zt^{n-1}$ with $u(1,x)\in V_n$.
Since $\chi_y\bigl((u(1,x)\bigr)$ does not depend on $y_1$
we can view $\chi_y\vert V_{n-1}$ as an element of $V_{n-1}^\ast$.
If we let $w_d^{n-1}\in V_{n-1}\subset V_n$
be the ${\Cal F}_{n-1}$-invariant elements as in Theorem
\thmref{\thmexplictcdescFinveltsinV}, then for $0\leq d\leq n-1$
$$
w_d^n= {{1-q}\over{1-q^n}} \sum_{l=0}^{n-1}
\r(T_lT_{l-1}\ldots T_2T_1) \, w^{n-1}_d
$$
by choosing minimal coset representatives in $S_n/S_{n-1}$.
For arbitrary $y\in\Ztn$ we get
$$
\phi_{w(y)}(w_d^n) = {{1-q}\over{1-q^n}} \sum_{l=0}^{n-1}
\bigl( \r^\ast(T_1T_2\ldots T_{l-1}T_l)\chi_y\bigr) (w_d^{n-1}).
$$
Using Theorem \thmref{\thmcontragrepone} we can calculate
$\r^\ast(T_1T_2\ldots T_{l-1}T_l)\chi_y$ explicitly for suitably
choosen $y\in\Ztn$, and we obtain
an explicit recurrence relation expressing $\phi_f(w_d^n)$
in terms of $\phi_f(w_d^{n-1})$ and $\phi_{f-1}(w_d^{n-1})$.
Together with the initial conditions for $\phi_0(w_d^n)$ and
$\phi_n(w_d^n)$ this recurrence relation determines
$\phi_f(w_d^n)$. This recurrence relation is equivalent to the
following contiguous relation for the $q$-Krawtchouk polynomials;
$$
\multline
(1-q^N)(1+aq^{N-2n})\, K_n(q^{-x};a,N;q) =
(1-q^{N-n})(1+aq^{N-n}) \, K_n(q^{-x};a,N-1;q)\\ +
q^{N-n}(1-q^n)(1+aq^{-n})\, K_{n-1}(q^{-x};a,N-1;q).
\endmultline
$$
\enddemo

In the next two remarks we discuss how Theorem
\thmref{\thmsphericalasqKrawtchouk} is related to
known interpretations of ($q$-)Krawtchouk polynomials
on finite (hyper)groups.

\demo{Remark \theoremname{\remrelatieDunkl}} In the specialisation
$q=p=1$ we have $V_n=\C[\Ztn]$ and  $K_f(d;{1\over 2},n)$ are
the spherical functions on $H_n$ with respect to subgroup $S_n$. This
result goes back to Vere-Jones in 1971 (also in the context of
statistics) and Delsarte in 1973 (related to the Hamming scheme
in coding theory), see \cite{\Dunk} and references given there.

In case we specialise only $q=1$, we see that $V_n$ is the $n$-fold
tensor product of $V_1$. Here $V_1$ is generated by two elements; a
unit element and $\omega$, say, satisfying $\omega^2 = (p-1)\omega+p$.
Then this can be considered as a commutative hypergroup, see Dunkl and
Ramirez \cite{\DunkR, \S 5} with their $a$ corresponding to $-p$.
The interpretation of Krawtchouk polynomials as symmetrised
characters \cite{\DunkR, Thm.~5.1} corresponds to the result in
Theorem \thmref{\thmsphericalasqKrawtchouk} specialised to $q=1$,
see also Koornwinder \cite{\Koor, \S 6}.
In this case the coefficients $c_l(k,d)$ occuring in the product
formula of Remark \thmref{\remproductformula} can be determined by
a counting argument and leads to the product formula
for Krawtchouk polynomials, see \cite{\DunkR, \S 5}.
In the general case $V_n$ gives rise to a $2^n$-point hypergroup,
which is not a $n$-fold tensor product. Its characters are described
in \S 5. 
\enddemo

\demo{Remark \theoremname{\remrelationStanton}}
Stanton \cite{\StanAJM}, \cite{\StanAKS}
has shown that the $q$-Krawtchouk polynomials appear as spherical
functions on certain finite groups of Lie type for specific
values of $p$ and $q$; see Carter \cite{\Cart} for information
on finite groups of Lie type. For a finite group $G$
with subgroup $B$ the Hecke algebra $H(G,B)$ is defined as
the algebra $e_B\C[G]e_B$, where $e_B = |B|^{-1}\sum_{b\in B}b$ is
the idempotent corresponding to the subgroup $B$. This can
also be viewed as the convolution algebra of the left and
right $B$-invariant functions on $G$, or as the
intertwiner algebra of the induced representation
$1_B^G$. In particular, $(G,B)$ is a Gelfand pair if and only
if $H(G,B)$ is commutative in which case we would like
to know the spherical functions. See \cite{\Cart}, \cite{\CurtRone},
\cite{\CurtRtwo} for more information and references.

Let $G$ be a finite group of Lie type, then it has a $BN$-pair
implying the existence of subgroups $B$ and $N$ such that there
is a Weyl group $W\cong N/(B\cap N)$. Moreover, the Bruhat
decomposition holds, $G= \cup_{w\in W} BwB$, and from this
we can associate to each parabolic subgroup of $W$ a
parabolic subgroup of $G$.
We now consider the cases in which the corresponding
Weyl group is the hyperoctahedral group $H_n$. Using
the classification of simple finite groups of Lie type
\cite{\Cart, \S 1.19, p.~464} there are $5$ types;
the Chevalley groups of type $B_n$ and $C_n$ and the
twisted groups of type ${}^2D_{n+1}$, ${}^2A_{2n-1}$ and
${}^2A_{2n}$. These groups have realisations as
classical groups over finite fields, cf. \cite{\Cart, p.~40}.
Moreover,
$H(G,B)$ is obtained from the generic Hecke algebra
${\Cal H}_n$ by suitable specialisation  of $p$ and $q$.
Denote by $P$ the corresponding
maximal parabolic subgroup corresponding to the
maximal parabolic subgroup $S_n\subset H_n$.
Then $H(G,P)\cong H(H_n,S_n)$, which is commutative,
see Curtis, Iwahori and
Kilmoyer \cite{\CurtIK, \S\S 2, 3} and also
\cite{\BrouCN, Thm.~10.4.11}.

The spherical functions corresponding to the
Gelfand pair $(G,P)$ have been determined by Stanton
\cite{\StanAJM} in terms of $q$-Krawtchouk polynomials
\thetag{\vgldefqKrawtchouk}; let $p_0(\not= 2)$ be a
prime and $q_0$ be an integral power of $p_0$, then the
zonal spherical functions can be expressed in terms of
$q$-Krawtchouk polynomials $K_f(\cdot ;p,n;q)$
with $p$ and $q$ as in the following table.
$$
\matrix
     & B_n & C_n & {}^2D_{n+1} & {}^2A_{2n-1} & {}^2A_{2n} \\
p    & q_0 & q_0 & q_0^2       & q_0          & q_0^3      \\
q    & q_0 & q_0 & q_0         & q_0^2        & q_0^2
\endmatrix
$$
(For these specialisations the values of the coefficients
$c_l(k,d)$ occuring in the product formula
of Remark \thmref{\remproductformula} can be determined,
since spherical functions satisfy a product formula,
see \cite{\StanAJM, \S 7}.)
These specialisations are also the specialisations of $p$ and $q$
needed to obtain $H(G,B)$ from ${\Cal H}_n$ in these five cases
according to \cite{\Cart, p.~464}, so
Theorem \thmref{\thmsphericalasqKrawtchouk}
unifies and generalises Stanton's results.

The reason for this is the following.
$V_n\cong V_n^\ast \cong e_B\C[G]e_P = \bigoplus_{f=0}^n U_f$
as ${\Cal H}_n=H(G,B)$ module and the action of the characteristic
function of $BwB$ can be considered as $\lambda(Bw)$ on
$L(B\backslash G/P)=e_B\C[G]e_P$, where $\lambda$ denotes the left
regular representation. By Stanton \cite{\StanAJM, Thm.~5.4},
which corresponds nicely to Theorem \thmref{\thmLaplacevgl}
and Corollary \thmref{\corthmLaplacevgl}, $L(G/P)=\C[G]e_P$
decomposes as $\bigoplus_{f=0}^n X_f$ under the left
regular representation $\lambda$ of $G$. Then $U_f = e_BX_f$
and $U_f^{{\Cal F}_n} = e_PU_f = e_PX_f$ is the one dimensional
space spanned by the zonal spherical function.
\enddemo

\demo{Remark \theoremname{\remrelationMacdonald}} Macdonald
\cite{\MacdSB} studies the representation of the (extended)
affine Hecke algebra obtained from inducing the index
representation of the corresponding finite
Hecke algebra for the Weyl group $W_0$.
The representation space can be identified with
the group algebra of the weight lattice $P$, and there exist
orthogonal polynomials $E_\lambda$ ($\lambda\in P$) acting on $P$
and Weyl group invariant orthogonal polynomials $P_\lambda$
($\lambda\in P^+$, the dominant weights), where the $P_\lambda$
can be obtained from Hecke symmetrising over $W_0$ from $E_\mu$
for $\mu$ in the Weyl group orbit of $\lambda$. The situation in
this paper is analogous to the situation
in Macdonald \cite{\MacdSB}; $P$,
$E_\lambda$, $P_\lambda$, $W$ correspond to
$\Ztn$, $\chi_y$, $\phi_f$, $S_n$.
In general, it seems not possible to derive the results of this
paper from Macdonald \cite{\MacdSB}. but see Matsumoto \cite{\Mats}
for the case $n=2$.
\enddemo

\demo{Remark \theoremname{\remadditionformulas}} As mentioned in
Remarks \thmref{\remrelatieDunkl} and
\thmref{\remrelationStanton}, Theorem
\thmref{\thmsphericalasqKrawtchouk} covers
a number of cases in which $\phi_f$ is a spherical function on
a finite group. In most of these cases this can be used to derive
an addition formula for ($q$-)Krawtchouk polynomials, see
Dunkl \cite{\Dunk} and Stanton \cite{\StanGD}. The interpretation
of general $q$-Krawtchouk polynomials as `zonal spherical functions'
on the Hecke algebra ${\Cal H}_n$ might lead to an addition formula
as well, but this is not clear.
It should be noted that there is also an addition formula for
so-called quantum $q$-Krawtchouk polynomials, which is derived
by Groza and Kachurik \cite{\GrozK}
from their relation to matrix elements of irreducible
representations of the quantum $SU(2)$ group and explict knowledge
of the Clebsch-Gordan coefficients for this quantum group.
This interpretation of quantum $q$-Krawtchouk polynomials on the
quantum $SU(2)$ group is different from the interpretation of
$q$-Krawtchouk polynomials on the quantum $SU(2)$ group as
described in Koornwinder \cite{\KoorZSE}, cf. Corollary
{\corthmLaplacevgl}.
\enddemo


\enddocument